\documentclass{amsart}

\usepackage{a4,xy,amssymb,diagrams,times}
\xyoption{poly}
\xyoption{arc}
\xyoption{knot}
\xyoption{curve}
\xyoption{arrow}
\xyoption{all}

\makeindex

\newcommand{\V}{\mathbb{V}}

\newcommand{\dummy[1]}{{#1}}
\newcommand{\C}{\mathbb{C}}

\newcommand{\R}{\mathbb{R}}
\newcommand{\N}{\mathbb{N}}
\newcommand{\Z}{\mathbb{Z}}

\newcommand{\PP}{\mathbb{P}}

\newcommand{\Ascr}{{\mathcal A}}
\newcommand{\Oscr}{{\mathcal O}}

\newcommand{\Af}{\mathbb{A}}

 \newcommand{\vtx}[1]{*+[o][F-]{\scriptscriptstyle #1}}

\newcommand{\wis}[1]{{\text{\em \usefont{OT1}{cmtt}{m}{n} #1}}}

\newtheorem{theorem}{Theorem}[section]
\newtheorem{lemma}[theorem]{Lemma}
\newtheorem{proposition}[theorem]{Proposition}

\theoremstyle{definition}
\newtheorem{definition}[theorem]{Definition}
\newtheorem{example}[theorem]{Example}

\theoremstyle{remark}

\title{Partial desingularizations arising from non-commutative algebras}

\author{Lieven Le Bruyn}
\address{Lieven Le Bruyn\\Universiteit Antwerpen\\B-2020 Antwerpen (Belgium)}
\email{lieven.lebruyn@ua.ac.be}
\author{Stijn Symens}
\address{Stijn Symens\\Universiteit Antwerpen\\B-2020 Antwerpen (Belgium)}
\email{stijn.symens@ua.ac.be}
\thanks{The second author is Research Assistant of the Fund for Scientific Research - Flanders (Belgium)}

\begin{document}

\sloppy
 
 \maketitle
 
 \begin{abstract}
 Let $X$ be a singular affine normal variety with coordinate ring $R$ and assume that there is an $R$-order $\Lambda$ admitting a stability structure $\theta$ such that the scheme of $\theta$-semistable representations is smooth, then we construct a partial desingularization of $X$ with classifiable remaining singularities. In dimension $3$ this explains the omnipresence of conifold singularities in partial desingularizations of quotient singularities. In higher dimensions we have a small list of singularity types generalizing the role of the conifold singularity.
 \end{abstract}
 
 \section{Introduction}
 
 In this paper we want to give a ringtheoretical explanation for the omnipresence of conifold singularities in partial desingularizations of three-dimensional quotient singularities coming from physics (see for example \cite{Berenstein} and \cite{BG98}) and to generalize this phenomenon to higher dimensions. For a translation between physics language and the mathematical terms used in this paper, we refer to section $4$ of our previous paper \cite{BLBS}.
 
 If $X=\C^3/G$ is a three-dimensional quotient singularity, one consider the McKay quiver setting $(Q,\alpha)$ of the finite group $G$ and the order over $\C[X]$
 \[
 \Lambda = \frac{\C Q}{R} \]
 obtained by dividing out commuting matrix relations, see for example \cite{CrawNotes}. One then chooses a stability structure $\theta$ such that the moduli space $\wis{moduli}^{\theta}_{\alpha}~\Lambda$ of isomorphism classes of $\theta$-semistable $\alpha$-dimensional $\Lambda$-representations is a partial resolution of $X$. In fact, in most examples, one even has that the scheme $\wis{rep}^{\theta-semist}_{\alpha}~\Lambda$ of $\theta$-semistable $\alpha$-dimensional representations is a smooth variety. In this paper we will show that this condition implies that possible remaining singularities in the (partial) desingularization
 \[
 \wis{moduli}^{\theta}_{\alpha}~\Lambda \rOnto X \]
 must be of conifold type. Moreover, we will extend this setting to higher dimensions.
 
 Let $X$ be an affine normal variety with coordinate ring $R = \C[X]$ and function field $K=\C(X)$. Let $\Lambda$ be an $R$-order in central simple $K$-algebra $\Sigma$ of dimension $n^2$. We say that $\Lambda$ is a {\em smooth $R$-order} if the scheme $\wis{trep}_n~\Lambda$ of trace preserving $n$-dimensional $\Lambda$-representations is a smooth variety. However, this is a very restrictive condition and usually an order $\Lambda$ will have a non-zero {\em defect} (to be defined in \S 2) to smoothness.
 
 Still, if $\Lambda$ has a complete set of orthogonal idempotents $\{ e_1,\hdots,e_k \}$ we have a well-defined dimension vector $\alpha=(a_1,\hdots,a_k) \in \N^k$ (where $a_i = tr_{\Lambda}(e_i)$) such that
 \[
 \wis{trep}_n~\Lambda \simeq \wis{GL}_n \times^{\wis{GL}(\alpha)} \wis{rep}_{\alpha}~\Lambda \]
 Let $\theta \in \Z^k$ such that $\theta.\alpha = 0$ then we define an $\alpha$-dimensional $\Lambda$-representation $V \in \wis{rep}_{\alpha}~\Lambda$ to be {\em $\theta$-semistable} if for all $\Lambda$-subrepresentations $W$ of $V$ we have $\theta.\beta \geq 0$ where $\beta$ is the dimension vector of $W$. The set of all $\alpha$-dimensional $\theta$-semistable representations $\wis{rep}_{\alpha}^{\theta-semist}~\Lambda$ is a Zariski open subset of $\wis{rep}_{\alpha}~\Lambda$.
 
 In favorable situations we can choose a stability structure $\theta$ such that $\wis{rep}_{\alpha}^{\theta-semist}~\Lambda$ is a smooth variety. In such a {\em good setting} we can use universal localization in $\wis{alg@n}$ to construct of sheaf $\mathcal{A}$ of smooth orders over the corresponding moduli space $\wis{moduli}_{\alpha}^{\theta}~\Lambda$ (parametrizing isomorphism classes of semistable $\alpha$-dimensional representations) giving a commutative diagram
 \[
\xymatrix@R=40pt@C=45pt{
\wis{spec}~\Ascr \ar[d]_c \ar[rd]^{\phi} \\
\wis{moduli}^{\theta}_{\alpha}~\Lambda \ar@{->>}[r]^{\pi} & X = \wis{spec}~R
}
\]
Here, $\wis{spec}~\mathcal{A}$ is a non-commutative variety obtained by gluing affine non-commutative structure sheaves $(\wis{spec}~\Lambda_D,\mathcal{O}^{nc}_{\Lambda_D})$ together. The map $c$ is defined locally by intersecting a prime ideal with its center and $\pi$ is a projective morphism. As $\mathcal{A}$ is a sheaf of smooth orders, one can view the resulting map $\phi$ as a {\em non-commutative desingularization} of $X$.

A good setting $(\Lambda,\alpha,\theta)$ also limits the types of remaining singularities in the partial desingularization $\pi$. If $\wis{dim}~X = 3$, the moduli space can have worst have conifold singularities, and in dimension $4,5$ resp. $6$ there is a full classification of the possible remaining singularities which consist of $4,10$ resp. $53$ types, see \cite{RBLBVdW}.

In the final section we study the special case of the conifold singularity in great detail. We give several ringtheoretical interpretations of the {\em conifold algebra} $\Lambda_c$ : as a skew-group ring over a polynomial ring and as a Clifford algebra. The latter description allows us to study the prime ideal structure of $\Lambda_c$ in great detail and determine its non-commutative structure sheaf $\mathcal{O}^{nc}_{\Lambda_c}$. We work out its scheme of $2$-dimensional representations, study the corresponding stability structures and work out the resulting desingularizations which are related by the so-called Atiyah flop.

The results contained in this paper were presented at the conference  'Sch\'emas de Hilbert, alg\`ebre non-commutative et correspondance de McKay' at CIRM, Luminy in october 2003, see \cite{LBnotes} for the lecture notes.

 \section{Geometry of orders}
 
 Let $X$ be a commutative normal variety  with affine coordinate ring the normal domain $R = \C[X]$ and function field $K = \C(X)$. Let $\Sigma$ be a central simple $K$-algebra of dimension $n^2$ and let $\Lambda$ be an $R$-order in $\Sigma$, that is, $\Lambda$ is an $R$-subalgebra of $\Sigma$ which is finitely generated as an $R$-module and such that $\Lambda.K = \Sigma$. Recall that there is a reduced trace map $tr : \Sigma \rTo K$ satisfying $tr(\Lambda) = R$ (because $R$ is integrally closed). Composing $tr$ with the inclusion $R \subset \Lambda$ we get a linear map $tr_{\Lambda} : \Lambda \rTo \Lambda$. In particular, if $\Lambda = M_n(R)$ the usual trace map induces the linear map $tr_{M_n(R)} : M_n(R) \rTo M_n(R)$ sending a matrix $A \in M_n(R)$ to the diagonal matrix $tr(A) 1_n$.
 
 The {\em scheme of trace preserving representations} $\wis{trep}_n~\Lambda$ is the affine scheme representing the functor
$\wis{commalg} \rTo \wis{sets}$ determined by
\[ \wis{trep}_n~\Lambda(\C) =  \{ \Lambda \rTo^{\phi} M_n(\C)~|~\text{$\phi$ an algebra morphism and~} \phi \circ tr_{\Lambda} = tr_{M_n(\C)} \circ \phi~\}.
\]
It is well known, see for example \cite{ProcesiCH} that conjugation of $M_n(\C)$ by $\wis{GL}_n(C)$ makes $\wis{trep}_n~\Lambda$ into an affine $\wis{GL}_n$-variety such that the corresponding algebraic quotient map
\[
\wis{trep}_n~\Lambda \rOnto \wis{trep}_n~\Lambda / \wis{GL}_n = \wis{triss}_n~\Lambda \simeq X = \wis{spec}~R \]
recovers the central variety $X$. One can also recover the order $\Lambda$ from the scheme of trace preserving representations as the algebra of $\wis{GL}_n$-equivariant maps from $\wis{trep}_n~\Lambda$ to $M_n(\C) = \Af^{n^2}_{\C}$ where the latter variety is a $\wis{GL}_n$-variety  under the action by conjugation, see again \cite{ProcesiCH}. The notation $\wis{triss}_n~\Lambda$ is motivated by the fact that the algebraic quotient of $\wis{trep}_n~\Lambda$ by $\wis{GL}_n$ classifies isomorphism classes of $n$-dimensional (trace preserving) semi-simple representations of $\Lambda$. That is, if $\mathfrak{m} \triangleleft R$ is a maximal ideal of $R$ with corresponding geometric point $x_{\mathfrak{m}} \in X$, then $\mathfrak{m}$ determines an $n$-dimensional semi-simple $\Lambda$-module
\[
M_{\mathfrak{m}} = S_1^{\oplus e_1} \oplus \hdots \oplus S_k^{\oplus e_k} ,
\]
where the $S_i$ are simple $\Lambda$-modules of dimension $d_i$ (and occurring in $M_{\mathfrak{m}}$ with multiplicity $e_i$) such that $\sum d_ie_i = n$. Indeed, the geometric point $x_{\mathfrak{m}}$ determines a trace preserving algebra map
\[
\overline{\Lambda}_{\mathfrak{m}} = \Lambda/ \mathfrak{m}\Lambda \rTo M_n(\C) \]
and hence an $n$-dimensional $\Lambda$-module $N_{\mathfrak{m}}$. The semi-simple module $M_{\mathfrak{m}}$ is the semi-simplification of $N_{\mathfrak{m}}$ that is the direct sum of its Jordan-H\"older factors. We say that $\mathfrak{m}$ (or the point $x_{\mathfrak{m}} \in X$) is of representation-type $\tau(\mathfrak{m}) = (e_1,d_1;\hdots;e_k,d_k)$.

To the maximal ideal $\mathfrak{m}$ we will associate a combinatorial tool, a quiver-setting $(Q_{\mathfrak{m}},\alpha_{\mathfrak{m}})$ where $Q_{\mathfrak{m}}$ is the quiver on $k$ vertices (with vertex $v_i$ corresponding to the simple component $S_i$) such that the number of oriented arrows from vertex $v_i$ to vertex $v_j$ is given by
\[
\#~\{ \xymatrix{\vtx{v_i} \ar[r] & \vtx{v_j}} \} = \wis{dim}_{\C}~Ext^1_{\Lambda}(S_i,S_j) \]
and where the dimension vector $\alpha_{\mathfrak{m}} = (e_1,\hdots,e_k)$ is determined by the multiplicities.  By this construction we have that the space of $\alpha_{\mathfrak{m}}$-dimensional representations of $Q_{\mathfrak{m}}$, $\wis{rep}_{\alpha_{\mathfrak{m}}} Q_{\mathfrak{m}}$ can be identified with the self-extension space $Ext^1_{\Lambda}(M_{\mathfrak{m}},M_{\mathfrak{m}})$. Observe that the action of the automorphism group $Aut_{\Lambda}(M_{\mathfrak{m}}) = \wis{GL}_{e_1} \times \hdots \times \wis{GL}_{e_k} = \wis{GL}(\alpha_{\mathfrak{m}})$ on the self-extensions $Ext^1_{\Lambda}(M_{\mathfrak{m}},M_{\mathfrak{m}})$ coincides with the action of $\wis{GL}(\alpha_{\mathfrak{m}})$ on $\wis{rep}_{\alpha_{\mathfrak{m}}} Q_{\mathfrak{m}}$ by base-change.
By definition of self-extensions every representation $V \in \wis{rep}_{\alpha_{\mathfrak{m}}} Q_{\mathfrak{m}}$ determines an algebra map
\[
\Lambda \rTo^{\phi_V} M_n(\C[\epsilon]), 
\]
where $\C[\epsilon] = \C[x]/(x^2)$ is the algebra of dual numbers. The $\wis{GL}(\alpha_{\mathfrak{m}})$-subspace of $\wis{rep}_{\alpha_{\mathfrak{m}}} Q_{\mathfrak{m}}$ consisting of all trace preserving extensions, that is such that $tr_{M_n(\C[\epsilon])} \circ \phi_V = \phi_V \circ tr_{\Lambda}$ can again be identified with the representation space of a {\em marked} quiver setting $\wis{rep}_{\alpha_{\mathfrak{m}}} Q^{\dagger}_{\mathfrak{m}}$ where $Q^{\dagger}_{\mathfrak{m}}$ is the same quiver as $Q_{\mathfrak{m}}$ except that certain loops may be removed and that some other loops may acquire a marking by which we mean that a representation of $Q^{\dagger}_{\mathfrak{m}}$ in a marked loop corresponds to a trace zero matrix, see \cite{LBetale} for more details. The whole point of this construction is that the normal space in $M_{\mathfrak{m}}$ to the closed orbit $\Oscr(M_{\mathfrak{m}})$ in the trace preserving representation space $\wis{trep}_n~\Lambda$
\[
\frac{T_{M_{\mathfrak{m}}}~\wis{trep}_n~\Lambda}{T_{M_{\mathfrak{m}}} \Oscr(M_{\mathfrak{m}})} = N_{M_{\mathfrak{m}}} \simeq \wis{rep}_{\alpha_{\mathfrak{m}}} Q^{\dagger}_{\mathfrak{m}}
\]
can be identified with the representation space of the marked quiver and that the automorphism is one as $\wis{GL}(\alpha_{\mathfrak{m}}) = Stab(M_{\mathfrak{m}})$ modules. This fact allows us to define a numerical {\em defect} measuring the failure of smoothness of $\wis{trep}_n~\Lambda$ over the point $x_{\mathfrak{m}}$.

\begin{definition} The {\em defect} $\wis{def}_{\mathfrak{m}}~\Lambda$ of the $R$-order $\Lambda$ in the maximal ideal $\mathfrak{m}$ is defined to be 
\[
\wis{def}_{\mathfrak{m}}~\Lambda = 1-\chi(\alpha_{\mathfrak{m}},\alpha_{\mathfrak{m}}) - \# \{ \text{marked loops in $Q^{\dagger}_{\mathfrak{m}}$} \} - \wis{dim}~X,
\]
where $\chi : \Z^k \times \Z^k \rTo \Z$ is the Euler form of the quiver obtained from $Q^{\dagger}_{\mathfrak{m}}$ by forgetting the markings, that is, the entry $(i,j)$ of the matrix defining $\chi$ is equal to $\delta_{ij} - \# \{ \xymatrix{\vtx{v_i} \ar[r] & \vtx{v_j}} \}$.
\end{definition}

\begin{proposition} With notations as above, $\wis{def}_{\mathfrak{m}}~\Lambda \geq 0$ and the following statements are equivalent
\begin{enumerate}
\item{$\wis{def}_{\mathfrak{m}}~\Lambda = 0 $.}
\item{$\wis{trep}_n~\Lambda$ is a smooth variety in all points lying over $x_{\mathfrak{m}}$.}
\end{enumerate}
\end{proposition}

\begin{proof} As $\Lambda$ is an $R$-order in an $n^2$-dimensional central simple  $K$-algebra $\Sigma$, there is a Zariski open subset $\wis{azu}_n~\Lambda$ of $X$ of points $x_{\mathfrak{m}}$ such that $\overline{\Lambda}_{\mathfrak{m}} \simeq M_n(\C)$ (the so called Azumaya locus of $\Lambda$). Over $\wis{azu}_n~\Lambda$ the algebraic quotient map $\wis{trep}_n~\Lambda \rOnto X$ is a principal $\wis{PGL}_n$-fiber whence generically the trace preserving representation scheme has dimension
\[
\wis{dim}~\wis{trep}_n~\Lambda = \wis{dim}~X + n^2 - 1.
\]
On the other hand, the dimension of the tangent space to the representation scheme in the semi-simple representation $M_{\mathfrak{m}}$ is equal to
\[
\begin{split}
\wis{dim}~T_{M_{\mathfrak{m}}}~\wis{trep}_n~\Lambda &= \wis{dim}~\Oscr(M_{\mathfrak{m}}) + \wis{dim}~\wis{rep}_{\alpha_{\mathfrak{m}}}~Q^{\dagger}_{\mathfrak{m}} \\
&= (n^2 - \sum_i e_i^2) + (\sum_{\xymatrix{\vtx{v_i} \ar[r] & \vtx{v_j}}} e_ie_j - \# \{ \text{marked loops in $Q^{\dagger}_{\mathfrak{m}}$} \} ) \\
&= n^2 - \chi(\alpha_{\mathfrak{m}},\alpha_{\mathfrak{m}}) - \# \{ \text{marked loops in $Q^{\dagger}_{\mathfrak{m}}$} \}
\end{split}
\]
and as $\wis{dim}~T_{M_{\mathfrak{m}}}~\wis{trep}_n~\Lambda \geq \wis{dim}~\wis{trep}_n~\Lambda$ it follows that $\wis{def}_{\mathfrak{m}}~\Lambda \geq 0$. Moreover, it also follows that $\wis{def}_{\mathfrak{m}}~\Lambda = 0$ if and only if $\wis{trep}_n~\Lambda$ is smooth in $M_{\mathfrak{m}}$. But as the singularities of $\wis{trep}_n~\Lambda$ form a $\wis{GL}_n$-closed subvariety and as $\Oscr(M_{\mathfrak{m}})$ is the unique closed orbit lying over $x_{\mathfrak{m}}$ (recall that closed orbits in $\wis{trep}_n~\Lambda$ are precisely the isomorphism classes of semi-simple representations) the equivalence of the two statements follows.
\end{proof}

\begin{example} Consider the quantum plane of order two $\Lambda = \C_{-1}[x,y]$ determined by the commutation relation $xy+yx=0$. If $u=x^2$ and $v=y^2$ then the center of $\Lambda$ is the polynomial algebra $R=\C[u,v]$ and $\Lambda$ is a free module of rank $4$ over it. In fact, $\Lambda$ is an $R$-order in the quaternion-algebra 
\[
\Sigma = \begin{pmatrix} u & & v \\ & \C(u,v) & \end{pmatrix} .
\]
The reduced trace map is determined by its images on a $\C$-basis
\[
tr(x^iy^j) = \begin{cases} 0 & \text{if either $i$ or $j$ is odd} \\
2x^iy^j & \text{if both $i$ and $j$ are even.}
\end{cases}
\]
In the affine plane $\Af^2 = \wis{spec}~R$ the Azumaya locus of $\Lambda$ is $\wis{azu}_2~\Lambda = \mathbb{X}(uv)$ the complement of the two coordinate axes. Let $x_{\mathfrak{m}} = (a^2,b) \in \mathbb{X}(uv)$ then the corresponding $2$-dimensional simple representation $M_{\mathfrak{m}}$ is determined by
\[
\Lambda \rOnto^{\phi} M_2(\C) \qquad \text{with} \qquad \phi(x) = \begin{bmatrix} a & 0 \\ 0 & -a \end{bmatrix} \qquad \phi(y) = \begin{bmatrix} 0 & 1 \\ b & 0 \end{bmatrix}.
\]
One verifies that $Ext^1_{\Lambda}(M_{\mathfrak{m}},M_{\mathfrak{m}}) \simeq \C^2$ and that the corresponding algebra map $\Lambda \rTo^{\psi}  M_2(\C[\epsilon])$ corresponding to $(\alpha,\beta) \in \C^2$ is given by
\[
\begin{cases}
\psi(u) &= \begin{bmatrix} a + \epsilon \alpha & 0 \\ 0 & -a - \epsilon \alpha \end{bmatrix} \\
\psi(v) &= \begin{bmatrix} 0 & 1 \\ b+\epsilon \beta & 0 \end{bmatrix}
\end{cases}
\]
and hence is trace preserving whence the local (marked) quiver-setting $(Q^{\dagger}_{\mathfrak{m}},\alpha_{\mathfrak{m}})$ is given by
\[
\xymatrix{\vtx{1} \ar@(l,ul) \ar@(ur,r)} \]
whence the defect is equal to $\wis{def}_{\mathfrak{m}}~\Lambda = 1 - (-1) - 0 - 2 = 0$ consistent with the fact that over the Azumaya locus (which is a smooth subvariety of the central scheme in this case) the algebraic quotient map is a principal $\wis{PGL}_2$-fibration whence $\wis{trep}_2~\Lambda$ will be smooth over it. For general orders $\Lambda$, if $x_{\mathfrak{m}}$ is a smooth point of the central variety and lies in the Azumaya locus, then $\wis{def}_{\mathfrak{m}}~\Lambda = 0$.

For $x_{\mathfrak{m}} = (a^2,0) \in \Af^2$ with $a \not= 0$ (and by a similar argument for points $(0,b)$ with $b \not= 0$), the corresponding semi-simple representation has two non-isomorphic one-dimensional simple components
\[
M_{\mathfrak{m}} = S_1 \oplus S_2 \qquad \text{with} \qquad S_i = \begin{cases} x \mapsto (-1)^i a \\ y \mapsto 0. \end{cases}
\]
One verifies that $Ext^1_{\Lambda}(S_i,S_i) = \C$ and that $Ext^1_{\Lambda}(S_1,S_2) \simeq Ext^1_{\Lambda}(S_2,S_1) \simeq \C$ whence the quiver-setting $(Q_{\mathfrak{m}},\alpha_{\mathfrak{m}})$ is given by
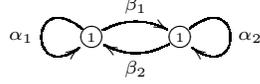
\[
\xymatrix{\vtx{1} \ar@(ul,dl)_{\alpha_1} \ar@/^/[r]^{\beta_1} & \vtx{1} \ar@(ur,dr)^{\alpha_2} \ar@/^/[l]^{\beta_2}} \]
and the corresponding algebra map $\Lambda \rTo M_2(\C[\epsilon])$ is given by
\[
x \mapsto \begin{bmatrix} a + \epsilon \alpha_1 & 0 \\ 0 & -a + \epsilon \alpha_2 \end{bmatrix} \qquad 
y \mapsto \begin{bmatrix} 0 & \beta_1 \\ \beta_2 & 0 \end{bmatrix} \]
which is only trace preserving if $\alpha_2 = - \alpha_1$ so we have one linear relation among the representations and therefore the corresponding (marked) quiver-setting $(Q^{\dagger}_{\mathfrak{m}},\alpha_{\mathfrak{m}})$ is equal to
\[
\xymatrix{\vtx{1} \ar@(ul,dl) \ar@/^/[r] & \vtx{1} \ar@/^/[l]} \]
and the defect is equal to $\wis{def}_{\mathfrak{m}}~\Lambda = 1 - (-1) - 0 - 2 = 0$ whence also over these ramified points the trace preserving representation variety $\wis{trep}_2~\Lambda$ is smooth.

Remains the point $x_{\mathfrak{m}} = (0,0)$ where the corresponding semi-simple representation is the zero-representation $M_{\mathfrak{m}} = S_0^{\oplus 2}$ where $S_0$ is determined by $x \mapsto 0$ and $y \mapsto 0$. One verifies that $Ext^1_{\Lambda}(S_0,S_0) \simeq \C^2$ whence the quiver-setting $(Q_{\mathfrak{m}},\alpha_{\mathfrak{m}})$ is equal to
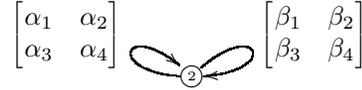
\[
\xymatrix{\vtx{2} \ar@(l,ul)^{\begin{bmatrix} \alpha_1 & \alpha_2 \\ \alpha_3 & \alpha_4 \end{bmatrix}} \ar@(ur,r)^{\begin{bmatrix} \beta_1 & \beta_2 \\ \beta_3 & \beta_4 \end{bmatrix}}}
\]
with corresponding algebra map $\Lambda \rTo M_2(\C[\epsilon])$ given by
\[
x \mapsto \epsilon \begin{bmatrix} \alpha_1 & \alpha_2 \\ \alpha_3 & \alpha_4 \end{bmatrix} \qquad y \mapsto \epsilon \begin{bmatrix} \beta_1 & \beta_2 \\ \beta_3 & \beta_4 \end{bmatrix} \]
which is only trace preserving if $\alpha_4 = - \alpha_1$ and $\beta_4 = - \beta_1$. Therefore the marked quiver-setting $(Q^{\dagger}_{\mathfrak{m}},\alpha_{\mathfrak{m}})$ is equal to
\[
\xymatrix{\vtx{2} \ar@(l,ul)|{\ast} \ar@(ur,r)|{\ast}} \]
and the defect is $\wis{def}_{\mathfrak{m}}~\Lambda = 1 -(-4)-2-2 = 1$ whence there must be a singularity of $\wis{trep}_2~\Lambda$ lying over $x_{\mathfrak{m}}$.

This is indeed  the case as the geometric points of $\wis{trep}_2~\Lambda$ are determined by couples of $2 \times 2 $ matrices
\[
( \begin{bmatrix} x_1 & x_2 \\ x_3 & -x_1 \end{bmatrix} , \begin{bmatrix} y_1 & y_2 \\ y_3 & -y_1 \end{bmatrix} ) \quad \text{satisfying} \quad tr( \begin{bmatrix} x_1 & x_2 \\ x_3 & -x_1 \end{bmatrix}.\begin{bmatrix} y_1 & y_2 \\ y_3 & -y_1 \end{bmatrix}) = 0. 
\]
That is, $\wis{trep}_2~\Lambda$ is the hypersurface in $\Af^6$ determined by the equation
\[
\wis{trep}_2~\Lambda = \mathbb{V}(2x_1y_1 + x_2y_3 + x_3y_2) \rInto \Af^6 \]
which is an irreducible $5$-dimensional variety having an isolated singularity at $x = (0,0,0,0,0,0)$ (the zero-representation).
\end{example}

\begin{definition} The {\em smooth locus} of an $R$-order $\Lambda$ is defined to be the subset of $X = \wis{spec}~R$
\[
\wis{smooth}_n~\Lambda = \{ x_{\mathfrak{m}} \in X~|~\wis{def}_{\mathfrak{m}}~\Lambda = 0 \}.
\]
We say that the order $\Lambda$ is {\em smooth} if $\wis{smooth}_n~\Lambda = X$, or equivalently, that $\wis{trep}_n~\Lambda$ is a smooth variety.
\end{definition}

If $X^{sm}$ denotes the smooth locus of $X = \wis{spec}~R$ then we have already seen that for any $R$-order $\Lambda$
\[
X^{sm} \cap \wis{azu}_n~\Lambda \rInto \wis{smooth}_n~\Lambda \]
as the algebraic quotient map $\wis{trep}_n~\Lambda \rOnto X$ is a principal $\wis{PGL}_n$-fibration over the Azumaya locus. In fact, for many interesting classes of orders the three loci coincide, that is,
\[
X^{sm} = \wis{azu}_n~\Lambda = \wis{smooth}_n~\Lambda.
\]
This is the case for quantum groups at roots of unity (see \cite{LBquantum}) and for orders associated at (deformed) preprojective algebras (see \cite{LBpreproj}). Later on we will prove a similar result for orders associated to quotient singularities.

If $x_{\mathfrak{m}} \in \wis{smooth}_n~\Lambda$ we know from \cite{LBetale} that the marked quiver setting $(Q^{\dagger}_{\mathfrak{m}},\alpha_{\mathfrak{m}})$ contains enough information to describe the \'etale local structure of $X$ near $x_{\mathfrak{m}}$ (that is, the structure of the $\mathfrak{m}$-adic completion $\hat{R}_{\mathfrak{m}}$) as well as the \'etale local structure of $\Lambda$ near $\mathfrak{m}$ (that is, the $\mathfrak{m}$-adic completion $\hat{\Lambda}_{\mathfrak{m}}$). We recall the result and refer to \cite{LBetale} for proof and more details.

\begin{proposition} Let $x_{\mathfrak{m}} \in \wis{smooth}_n~\Lambda$ with associated marked quiver-setting $(Q^{\dagger}_{\mathfrak{m}},\alpha_{\mathfrak{m}})$ with $\alpha_{\mathfrak{m}} = (a_1,\hdots,a_k)$. Then,
\begin{enumerate}
\item{The $\mathfrak{m}$-adic completion of the center $\hat{R}_{\mathfrak{m}}$ is isomorphic to the completion of the algebra generated by traces along oriented cycles in $(Q^{\dagger}_{\mathfrak{m}},\alpha_{\mathfrak{m}})$ at the maximal ideal generated by these traces.}
\item{The $\mathfrak{m}$-adic completion of the order $\Lambda$ is of the form
\[
\hat{\Lambda}_{\mathfrak{m}} \simeq \begin{bmatrix} M_{11} & \hdots & M_{1k} \\
\vdots & & \vdots \\
M_{k1} & \hdots & M_{kk} \end{bmatrix} \]
where $M_{ij}$ is a block of size $a_i \times a_j$ with all entries equal to the $\hat{R}_{\mathfrak{m}}$-module generated by all paths in $(Q^{\dagger}_{\mathfrak{m}},\alpha_{\mathfrak{m}})$ starting at vertex $v_i$ and ending in vertex $v_j$.}
\end{enumerate}
\end{proposition}

In particular, if $x_{\mathfrak{m}} \in \wis{smooth}_n~\Lambda$ we can describe the finite dimensional algebra $\overline{\Lambda}_{\mathfrak{m}} = \Lambda / \mathfrak{m} \Lambda$ to be Morita equivalent to the quotient of the path algebra of the underlying quiver $\C Q^{\dagger}_{\mathfrak{m}}$ by the ideal generated by all cycles in $Q^{\dagger}_{\mathfrak{m}}$.

\begin{definition}
Let $\wis{cat}$ be a category of $\C$-algebras. We say that an algebra $A \in \wis{cat}$ is {\em $\wis{cat}$-smooth} if and only if for every $B \in \wis{cat}$, every quotient $B \rOnto^{\pi} B/I$ in $\wis{cat}$ with $I$ a nilpotent ideal and every algebra morphism $A \rTo^{\phi} B/I$ in $\wis{cat}$ the diagram
\[
\xymatrix@R=45pt@C=45pt{
A \ar[rd]_{\phi} \ar@{.>}[r]^{\tilde{\phi}} & B \ar@{->>}[d]^{\pi} \\
& B/I 
}
\]
can be completed by an algebra morphism $A \rTo^{\tilde{\phi}} B$ in $\wis{cat}$.
\end{definition}

Grothendieck proved that an affine commutative $\C$-algebra $R$ is $\wis{commalg}$-smooth if and only if $R$ is regular, that is, if and only if $X = \wis{spec}~R$ is a smooth variety. Cuntz and Quillen \cite{CuntzQuillen} introduced {\em quasi-free algebras} as coordinate rings of non-commutative algebraic manifolds and they are precisely the $\wis{alg}$-smooth algebras. Similarly, smooth orders are $\wis{alg@n}$-smooth algebras where $\wis{alg@n}$ is the category of Cayley-Hamilton algebras of degree $n$ which we will describe briefly and refer to \cite{ProcesiCH} for more details.

If $M \in M_n(R)$ for $R$ a commutative $\C$-algebra, then its characteristic polynomial
\[
\chi_M = det(t1_n-M) = t^n + a_1 t^{n-1} + \hdots + a_n \]
is such that all its coefficients are polynomials with rational coefficients in traces of powers of $M$, that is, $a_i = f_i(Tr(M),Tr(M^2),\hdots,Tr(M^{n-1}))$. Hence, if $A$ is a $\C$-algebra having a trace map $tr_A~:~A \rTo A$ (a linear map satisfying $tr_A(tr_A(a)b)=tr_A(a)tr_A(b)$, $tr_A(ab)=tr_A(ba)$ and $tr_A(a)~b=b~tr_A(a)$ for all $a,b \in A$) then we define a {\em formal characteristic polynomial of degree $n$} for every $a \in A$ by
\[
\chi_a = t^n + f_1(tr_A(a),\hdots,tr_A(a^{n-1})) t^{n-1} + \hdots + f_n(tr_A(a),\hdots,tr_A(a^{n-1})) \]

\begin{definition} An object of $\wis{alg@n}$ is a Cayley-Hamilton algebra of degree $n$, that is, a $\C$-algebra having a trace map $tr_A$ satisfying
\[
\forall a \in A:~\chi_a(a) = 0~ \qquad \text{and} \qquad tr_A(1) = n \]
Morphisms $A \rTo^f B$ in $\wis{alg@n}$ are $\C$-algebra morphisms preserving traces, that is
\[
\xymatrix@R=45pt@C=45pt{
A  \ar[r]^f \ar[d]_{tr_A} & B \ar[d]_{tr_B} \\
A  \ar[r]^f & B
}
\]
is a commutative diagram.
\end{definition}

We recall from \cite{ProcesiCH} that $A \in \wis{alg@n}$ is $\wis{alg@n}$-smooth if and only if $\wis{trep}_n~A$ is a smooth variety (possibly having several irreducible components). In particular, a smooth order $\Lambda$ in a central simple $K$-algebra $\Sigma$ of dimension $n^2$ equipped with the reduced trace map is $\wis{alg@n}$-smooth.

Having identified smooth orders as a natural generalization of regular commutative algebras to the category of Cayley-Hamilton algebras and having a combinatorial local description of them (as well as their centers), we now turn to the associated {\em non-commutative smooth variety}.

\begin{definition}
Let $\Lambda$ be an $R$-order in a central simple $K$-algebra $\Sigma$ of dimension $n^2$, then the {\em non-commutative spectrum}, $\wis{spec}~\Lambda$ is the set of all twosided prime ideals $P$ of $\Lambda$ (that is, the ideals satisfying $a \Lambda b \subset P \Rightarrow a$ or $b \in P$). This set is equipped with the {\em Zariski topology} with typical open sets
\[
\mathbb{X}(I) = \{ P \in \wis{spec}~\Lambda~|~I \not\subset P \}
\]
for any twosided ideal $I$ of $\Lambda$
(see for example \cite{FVO444} and \cite{FVOAV}). The topological space $\wis{spec}~\Lambda$ comes equipped with a {\em non-commutative structure sheaf} $\mathcal{O}^{nc}_{\Lambda}$ with sections on the open set $\mathbb{X}(I)$
\[
\Gamma(\mathbb{X}(I),\mathcal{O}^{nc}_{\Lambda}) = \{ \delta \in \Sigma~|~\exists l \in \N~:~I^l.\delta \subset \Lambda \} \]
(again see \cite{FVO444} or \cite{FVOAV} for a proof that this defines a sheaf of non-commutative algebras with global sections $\Gamma(\wis{spec}~\Lambda,\mathcal{O}^{nc}_{\Lambda}) = \Lambda$).
Moreover, the {\em stalk} of $\mathcal{O}^{nc}_{\Lambda}$ at a prime ideal $P \in \wis{spec}~\Lambda$ is the symmetric localization
\[
\mathcal{O}^{nc}_{\Lambda,P} = Q_{\Lambda-P}(\Lambda) = \{ \delta \in \Sigma~|~I \delta \subset \Lambda~\text{for some twosided ideal}~I \not\subset P \}.
\]
\end{definition}

Intersecting a twosided prime ideal $P$ of $\Lambda$ with its center gives a prime ideal of $R$ and hence we obtain a continuous map
\[
\wis{spec}~\Lambda \rTo^{\pi_c} \wis{spec}~R \qquad P \mapsto P \cap R \]
and if we denote with $\Oscr_{\Lambda}$ the (usual) sheaf of $R$-algebras on $\wis{spec}~R$ associated to the $R$-order $\Lambda$ then $\pi_c$ induces a morphism of sheaves of algebras
\[
(\wis{spec}~\Lambda, \Oscr^{nc}_{\Lambda}) \rTo^{\pi_c} (\wis{spec}~R, \Oscr_{\Lambda}).
\]
For $\mathfrak{m}$ a maximal ideal of $R$ we can relate the local marked quiver setting $(Q^{\dagger}_{\mathfrak{m}},\alpha_{\mathfrak{m}})$ to the fiber $\pi_c^{-1}(\mathfrak{m})$. This quiver setting was determined by the semi-simple $n$-dimensional $\Lambda$-representation
\[
M_{\mathfrak{m}} = S_1^{\oplus e_1} \oplus \hdots \oplus S_k^{\oplus e_k} \]
where $S_i$ is a simple $d_i$-dimensional $\Lambda$-representation. Then, we have that
\[
\pi_c^{-1}(\mathfrak{m}) = \{ P_1,\hdots,P_k \} \qquad \text{with} \qquad \Lambda/P_i \simeq M_{d_i}(\C) \]
so the number of vertices in $Q^{\dagger}_{\mathfrak{m}}$ determines the number of maximal twosided ideals of $\Lambda$ lying over $\mathfrak{m}$ and the dimension vector $\alpha_{\mathfrak{m}} = (e_1,\hdots,e_k)$ determines the so called Bergman-Small data, see \cite{BergmanSmall}. The finitely many maximal twosided ideals $\{ P_1,\hdots,P_k \}$ lying over the central point $\mathfrak{m}$ form a {\em clique} \cite{Jategaonkar} and should be thought of as points lying infinitesimally close together in $\wis{spec}~\Lambda$. The marked quiver $Q^{\dagger}$ encodes this infinitesimal information. If $\mathfrak{m}$ is a central singularity, the hope is that one can use these finitely many infinitesimally close points to separate tangent information in $\mathfrak{m}$ rather than having to resort to the full blown-up of $\mathfrak{m}$. In the next section we will give some examples when this non-commutative approach to desingularization actually works.

\begin{example} Let $X = \Af^1$, that is $R=\C[x]$ and consider the order
\[
\Lambda = \begin{bmatrix} R & R \\ \mathfrak{m} & R \end{bmatrix},
\]
where $\mathfrak{m} = (x) \triangleleft R$, that is $x_{\mathfrak{m}} = 0$. For every point $\lambda \not= 0$ there is a unique maximal twosided ideal of $\Lambda$ lying over $\mathfrak{m}_{\lambda} = (x-\lambda)$ with quotient $M_2(\C)$. For this reason we say that $X - \{ 0 \}$ is the {\em Azumaya locus} of $\Lambda$. On the other hand, the {\em ramification locus} of $\Lambda$ is the closed subset $\{ 0 \} = \mathbb{V}(x)$ and there are two maximal ideals of $\Lambda$ lying over $\mathfrak{m}$
\[
M_1 = \begin{bmatrix} \mathfrak{m} & R \\ \mathfrak{m} & R \end{bmatrix} \qquad \text{and} \qquad M_2 = \begin{bmatrix} R & R \\ \mathfrak{m} & \mathfrak{m} \end{bmatrix} \]
and the quotients are $\Lambda/M_1 \simeq \C \simeq \Lambda/M_2$ whence they determine both a one-dimensional $\Lambda$-representation. That is, the canonical continuous map
\[
\wis{spec}~\Lambda \rOnto^{\pi_c} \wis{spec}~R \]
is a homeomorphism over $\mathbb{X}(x)$ and there are precisely two (infinitesimally close) points lying over $\mathbb{V}(x)$. The corresponding (marked) quiver setting is
\[
\xymatrix{\vtx{1} \ar@/^/[r] & \vtx{1} \ar@/^/[l]}
\]
and so the defect $\wis{def}_{\mathfrak{m}}~\Lambda = 0$. Remark that in all other maximal ideals $\mathfrak{m}_{\lambda}$ the local (marked) quiver setting is
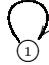
\[
\xymatrix{\vtx{1} \ar@(ul,ur)} \]
which also has zero defect so $\Lambda$ is a smooth order and hence $\wis{trep}_2~\Lambda$ is a smooth variety. We now turn to the structure sheaves $\Oscr_{\Lambda}$ and $\Oscr_{\Lambda}^{(nc)}$.
The central structure sheaf is just given by central localization and therefore we find for its stalks
\[
\Oscr_{\Lambda,\mathfrak{m}} = \begin{bmatrix} R_{\mathfrak{m}} & R_{\mathfrak{m}} \\ R_{\mathfrak{m}} & R_{\mathfrak{m}} \end{bmatrix} \qquad \Oscr_{\Lambda,\mathfrak{m}_{\lambda}} \simeq \begin{bmatrix} R_{\mathfrak{m}_{\lambda}} & R_{\mathfrak{m}_{\lambda}} \\ R_{\mathfrak{m}_{\lambda}} & R_{\mathfrak{m}_{\lambda}} \end{bmatrix}.
\]
Over the Azumaya locus the non-commutative structure sheaf $\Oscr_{\Lambda}^{nc}$ coincides with the central structure sheaf. The stalks in the two points lying over $\mathfrak{m}$ can be computed to be
\[
\Oscr_{\Lambda,M_1}^{nc} \simeq \begin{bmatrix} R_{\mathfrak{m}} & R_{\mathfrak{m}} \\ R_{\mathfrak{m}} & R_{\mathfrak{m}} \end{bmatrix} \qquad \Oscr_{\Lambda,M_2}^{nc} \simeq
\begin{bmatrix} R_{\mathfrak{m}} & x^{-1} R_{\mathfrak{m}} \\ x R_{\mathfrak{m}} & R_{\mathfrak{m}} \end{bmatrix}, \]
both of them being Azumaya algebras. Hence, we have the slightly surprising fact that the non-commutative structure sheaf $\Oscr_{\Lambda}^{nc}$ over $\wis{spec}~\Lambda$ is a sheaf of Azumaya algebras whereas $\Lambda$ itself is ramified in $\mathfrak{m}$. Observe that the stalk in $\mathfrak{m}$ of the central structure sheaf is the intersection of the two Azumaya stalks of the non-commutative structure sheaf.
\end{example}

\section{Moduli spaces}

In this section, $\Lambda$ will be an $R$-order in a central simple $K$-algebra of dimension $n^2$ and $\mathfrak{m}$ will be a singularity of $\wis{spec}~R = X$. W want to use $\Lambda$ to resolve the singularity in $\mathfrak{m}$. As we are only interested in the \'etale local structure of the singularity in $\mathfrak{m}$ we may restrict attention to $\hat{\Lambda}_{\mathfrak{m}}$ or more generally it is only the \'etale local structure of $\Lambda$ that is important. Hence, we may assume that $\Lambda$ is split as far as possible, or equivalently, that we have a complete set $\{ e_1,\hdots,e_k \}$ of orthogonal idempotents in $\Lambda$. That is the $e_i$ satisfy
\[
e_i^2 = e_i \qquad e_i.e_j = 0~\text{for $i \not= j$} \qquad \sum_{i=1}^k e_i = 1_{\Lambda} \]
These idempotents allow us to decompose finite dimensional $\Lambda$-representations. If $V \in \wis{rep}_m~\Lambda$ is an $m$-dimensional representation, we say that $V$ is of {\em dimension vector} $\alpha = (a_1,\hdots,a_n)$ for $\sum_{i=1}^k a_i = m$ provided
\[
\wis{dim}_{\C}~e_i.V = a_i \]
We denote this by $\wis{dim}~V = \alpha$.
Because $S = \overbrace{\C \times \hdots \times \C}^k \rInto \Lambda$ we can restrict $m$-dimensional $\Lambda$ representations to the semi-simple subalgebra $S$ to obtain morphisms
\[
\wis{rep}_m~\Lambda \rTo \wis{rep}_m~S = \bigsqcup_{\alpha}~\wis{GL}_m/\wis{GL}(\alpha) \]
where the decomposition is taken over all dimension vectors $\alpha = (a_1,\hdots.a_k)$ such that $\sum_i a_i = m$ and where $\wis{GL}(\alpha) = \wis{GL}_{a_1} \times \hdots \times \wis{GL}_{a_k}$. The component $\wis{GL}_m/\wis{GL}(\alpha)$ is the orbit of the semi-simple $S$-representation $V_{\alpha}$ with action given by the matrices
\[
e_i \mapsto E_{\sum_{j=1}^{i-1} a_j +1,\sum_{j=1}^{i-1} a_j + 1} + E_{\sum_{j=1}^{i-1} a_j +2,\sum_{j=1}^{i-1} a_j + 2} + \hdots + E_{\sum_{j=1}^{i} a_j,\sum_{j=1}^{i} a_j} \]
where $E_{i,j}$ are the standard matrices $(\delta_{iu}\delta_{jv})_{u,v} \in M_m(\C)$. As a consequence we can also decompose the representation schemes
\[
\wis{rep}_m~\Lambda = \bigsqcup_{\alpha} \wis{GL}_m \times^{\wis{GL}(\alpha)} \wis{rep}_{\alpha}~\Lambda \]
where $\wis{rep}_{\alpha}~\Lambda$ is the scheme representing all $m = \sum_i a_i$-dimensional representations of dimension vector $\alpha = (a_1,\hdots,a_k)$ on which the action by the set of idempotents $\{ e_1,\hdots,e_k \}$ is given by the above matrices. Clearly, the reductive group $\wis{GL}(\alpha)$ acts by base-change in the subspaces $e_i.V$ on $\wis{rep}_{\alpha}~\Lambda$ and the corresponding component of $\wis{rep}_m~\Lambda$ is the principal fiber bundle $\wis{GL}_m \times^{\wis{GL}(\alpha)} \wis{rep}_{\alpha}~\Lambda$.

A {\em character} of the reductive group $\wis{GL}(\alpha)$ is determined by an integral $k$-tuple $\theta = (t_1,\hdots,t_k) \in \Z^k$
\[
\chi_{\theta}~:~\wis{GL}(\alpha) \rTo \C^* \qquad (g_1,\hdots,g_k) \mapsto det(g_1)^{t_1} \hdots det(g_k)^{t_k} \]
As the subgroup $\C^*(1_{a_1},\hdots,1_{a_k})$ acts trivially on $\wis{rep}_{\alpha}~\Lambda$ we are only interested in the characters $\chi_{\theta}$ such that $0 = \theta.\alpha = \sum_{i=1}^k a_it_i$. Remark that a $\Lambda$-subrepresentation $W \subset V$ for $V \in \wis{rep}_{\alpha}~\Lambda$ necessarily satisfies $W \in \wis{rep}_{\beta}~\Lambda$ for some dimension vector $\beta \leq \alpha$.
We will now extend the definition of (semi)stable representations of quivers, due to A. King \cite{King} to the present setting.

\begin{definition} For $\theta \in \Z^k$ satisfying $\theta.\alpha = 0$, a representation $V \in \wis{rep}_{\alpha}~\Lambda$ is said to be
\begin{enumerate}
\item{{\em $\theta$-semistable} if and only if for every proper $\Lambda$-subrepresentation $W \subset V$ we have $\theta.\wis{dim}~W \geq 0$.}
\item{{\em $\theta$-stable} if and only if for every proper $\Lambda$-subrepresentation $W \subset V$ we have $\theta.\wis{dim}~W > 0$.}
\end{enumerate}
\end{definition}

For any setting satisfying $\theta.\alpha = 0$ we have the following inclusions of Zariski open $\wis{GL}(\alpha)$-stable subschemes of $\wis{rep}_{\alpha}~\Lambda$ (with obvious notations)
\[
\wis{rep}_{\alpha}^{simple}~\Lambda \subset \wis{rep}^{\theta-stable}_{\alpha}~\Lambda \subset \wis{rep}_{\alpha}^{\theta-semist}~\Lambda \subset \wis{rep}_{\alpha}~\Lambda \]
but some of these open subsets may actually be empty.

All these definitions carry over to any affine $\C$-algebra $\Lambda$ but if $\Lambda$ is an $R$-order in a central simple $K$-algebra of dimension $n^2$ we have the following link with the material of the previous section
\[
\wis{trep}_n~\Lambda = \wis{GL}_n \times^{\wis{GL}(\alpha)} \wis{rep}_{\alpha}~\Lambda \]
for the dimension vector $\alpha = (tr_{\Lambda}(e_1),\hdots,tr_{\Lambda}(e_k))$. Moreover,
\[
R = \C[\wis{triss}_n~\Lambda] = \C[\wis{iss}_{\alpha}~\Lambda] = \C[\wis{rep}_{\alpha}~\Lambda]^{\wis{GL}(\alpha)} \]
where $\wis{iss}_{\alpha}~\Lambda$ is the scheme representing semi-simple $\alpha$-dimensional representations of $\Lambda$. Remark that the dimension vector $\alpha$ above is such that there are $\alpha$-dimensional simple representations of $\Lambda$ so that in the above inclusion of $\wis{GL}(\alpha)$-stable subvarieties of $\wis{rep}_{\alpha}~\Lambda$ none of the subschemes is empty. From now on we fix this particular dimension vector $\alpha$ of total dimension $n$.

A polynomial function $f \in \C[\wis{rep}_{\alpha}~\Lambda]$ is said to be a {\em $\theta$-semi-invariant of weight $l$} if and only if we have for all $g \in \wis{GL}(\alpha)$
\[
g.f = \chi_{\theta}(g)^l f \]
where, as before, $\chi_{\theta}$ is the character of $\wis{GL}(\alpha)$ corresponding to $\theta$. It follows from \cite{King} that a representation $V \in \wis{rep}_{\alpha}~\Lambda$ is $\theta$-semistable if and only if there is some $\theta$-semi-invariant $f$ of some weight $l$ such that $f(V) \not= 0$.

Clearly, $\theta$-semi-invariants of weight zero are just polynomial invariants in $\C[\wis{rep}_{\alpha}~\Lambda]^{\wis{GL}(\alpha)} = R$ and the multiplication of $\theta$-semi-invariants of weights $l$ resp. $l'$ is a $\theta$-semi-invariant of weight $l+l'$. Therefore, the ring of all $\theta$-semi-invariants
\[
\C[\wis{rep}_{\alpha}~\Lambda]^{\wis{GL}(\alpha),\theta} = \bigoplus_{l=0}^{\infty} \{ f \in \C[\wis{rep}_{\alpha}~\Lambda]~|~\forall g \in \wis{GL}(\alpha)~:~g.f = \chi_{\theta}^l f \} \]
is a graded algebra with part of degree zero $R = \C[\wis{iss}_{\alpha}~\Lambda]$. Consequently, we have a projective morphism
\[
\wis{proj}~\C[\wis{rep}_{\alpha}~\Lambda]^{\wis{GL}(\alpha),\theta} \rOnto^{\pi} X = \wis{spec}~R \]
such that all fibers of $\pi$ are projective varieties. The main results of $\pi$ are proved as in \cite{King}.

\begin{theorem} There is a one-to-one correspondence between
\begin{enumerate}
\item{points in $\wis{proj}~\C[\wis{rep}_{\alpha}~\Lambda]^{\wis{GL}(\alpha),\theta}$, and}
\item{isomorphism classes of direct sums of $\theta$-stable $\Lambda$ representations of total dimension $\alpha$.}
\end{enumerate}
Moreover, as there are simple $\alpha$-dimensional $\Lambda$-representations, the morphism $\pi$ is  a birational projective map.
\end{theorem}

\begin{definition} We call $\wis{proj}~\C[\wis{rep}_{\alpha}~\Lambda]^{\wis{GL}(\alpha),\theta}$ the {\em moduli space of $\theta$-semistable representations} of $\Lambda$ and will denote it with $\wis{moduli}^{\theta}_{\alpha}~\Lambda$.
\end{definition}

Let us recall some examples of current interest.

\begin{example}[Kleinian singularities]
For a Kleinian singularity, that is, a quotient singularity  $\C^2/G$ with $G \subset SL_2(\C)$ there is an extended Dynkin diagram $D$ associated. 
Let $Q$ be the {\em double quiver} of $D$, that is to each arrow $\xymatrix{\vtx{} \ar[r]^x & \vtx{}}$ in $D$ we adjoin an arrow $\xymatrix{\vtx{} & \vtx{} \ar[l]_{x^*} }$ in $Q$ in the opposite direction and let $\alpha$ be the unique  minimal dimension vector such that $\chi_D(\alpha,\alpha) = 0$ (the so called isotropic Schur root of the tame quiver $\vec{D}$ obtained from the graph $D$ by fixing a certain orientation on the edges). Consider the {\em moment element}
\[
m = \sum_{x \in D} [x,x^*] \]
then the skew-group algebra $\Lambda = \C[x,y] \# G$ is on $R$-order with $R = \C[\C^2/G]$ in $M_n(K)$ where $K$ is the field of fractions of $R$ and $n = \#G$. Moreover, $\Lambda$ is Morita equivalent to the {\em preprojective algebra} which is the quotient of the path algebra of $Q$ by the ideal generated by the moment element
\[
\Pi_0 = \C Q/ (\sum [x,x^*] ).
 \]
For more details we refer to the lecture notes by W. Crawley-Boevey \cite{CrawleyLectNotes}. 
 If we take $\theta$ to be a generic character such that $\theta.\alpha = 0$, then the projective map
\[
\wis{moduli}^{\theta}_{\alpha}~\Lambda \rOnto X = \C^2/G \]
is a minimal resolution of singularities. Note that the map is birational as $\alpha$ is the dimension vector of a simple representation of $A = \Pi_0$, see \cite{CrawleyLectNotes}.

For such a stability structure $\theta$ we have that $\wis{rep}^{\theta-semist}_{\alpha}~\Pi_0$ is a smooth variety.
For consider the {\em moment map}
\[
\wis{rep}_{\alpha}~Q \rTo^{\mu} \wis{lie}~\wis{GL}(\alpha) = M_{\alpha}(\C) = M_{e_1}(\C) \oplus \hdots \oplus M_{e_k}(\C)
\]
defined by sending $V = (V_a,V_{a^*})$ to
\[
 (\sum_{\xymatrix{\vtx{} \ar[r]^a&\vtx{1}}} V_aV_{a^*} - \sum_{\xymatrix{\vtx{1} \ar[r]^a & \vtx{}}} V_{a^*}V_a, \hdots, \sum_{\xymatrix{\vtx{} \ar[r]^a&\vtx{k}} }V_aV_{a^*} - \sum_{\xymatrix{\vtx{k} \ar[r]^a & \vtx{}} }V_{a^*}V_a).\]
The differential $d \mu$ can be verified to be surjective in any representation $V \in \wis{rep}_{\alpha}~Q$ which has stabilizer subgroup $\C^*(1_{e_1},\hdots,1_{e_k})$ (a so called {\em Schur representation}) see for example \cite[lemma 6.5]{CrawleyMoment}. 

Further, any $\theta$-stable representation is Schurian. Moreover, for a generic stability structure $\theta \in \Z^k$ we have that every $\theta$-semistable $\alpha$-dimensional representation is $\theta$-stable as the $gcd(\alpha) = 1$.
Combining these facts it follows that $\mu^{-1}(0) = \wis{rep}_{\alpha}~\Pi_0$ is smooth in all $\theta$-stable representations.
\end{example}

\begin{example} 
Consider a quotient singularity $X = \C^d/G$ with $G \subset SL_d(\C)$ and $Q$ be the {\em McKay quiver} of $G$ acting on $V=\C^d$. 
That is, the vertices $\{ v_1,\hdots,v_k \}$ of $Q$ are in one-to-one correspondence with the irreducible representations $\{ R_1,\hdots,R_k \}$ of $G$ such that  $R_1 = \C_{triv}$ is the trivial representation. Decompose the tensorproduct in irreducibles
\[
V \otimes_{\C} R_j = R_1^{\oplus j_1} \oplus \hdots \oplus R_k^{\oplus j_k}, \]
then the number of arrows in $Q$ from $v_i$ to $v_j$
\[
\#~(v_i \rTo v_j ) = j_i \]
is the multiplicity of $R_i$ in $V \otimes R_j$. Let $\alpha = (e_1,\hdots,e_k)$ be the dimension vector where $e_i = \wis{dim}_{\C}~R_i$. 

The relevance of this quiver-setting is that
\[
\wis{rep}_{\alpha}~Q = Hom_G(R,R \otimes V) \]
where $R$ is the {\em regular representation}, see for example \cite{CrawNotes}. Consider $Y \subset \wis{rep}_{\alpha}~Q$ the affine subvariety of all $\alpha$-dimensional representations of $Q$ for which the corresponding $G$-equivariant map $B \in Hom_G(R,V \otimes R)$ satisfies
\[
B \wedge B = 0 \in Hom_G(R,\wedge^2 V \otimes R). \]
$Y$ is called the {\em variety of commuting matrices} and its defining relations can be expressed as linear equations between paths in $Q$ evaluated in $\wis{rep}_{\alpha}~Q$, say $(l_1,\hdots,l_z)$. Then, the quiver-order
\[
\Lambda = \frac{\int_{\alpha} \C Q}{(l_1,\hdots,l_z)} \]
is an order with center $R = \C[\C^d/G]$. In fact, $\Lambda$ is just the skew group algebra
\[
A = \C[x_1,\hdots,x_d] \# G. \]
Assume that the first vertex in the McKay quiver corresponds to the trivial representation. Take a character $\theta \in \Z^k$ such that $t_1 < 0$ and all $t_i > 0$ for $i \geq 2$, for example take
\[
\theta = ( - \sum_{i=2}^k \wis{dim} R_i , 1, \hdots, 1 ). \]
Then, the corresponding moduli space is isomorphic to
\[
\wis{moduli}^{\theta}_{\alpha}~A \simeq G-\wis{Hilb}~\C^d \]
the {\em $G$-equivariant Hilbert scheme} which classifies all $\# G$-codimensional ideals $I \triangleleft \C[x_1,\hdots,x_d]$ where
\[
\frac{\C[x_1,\hdots,x_d]}{I} \simeq \C G \]
as $G$-modules, hence in particular $I$ must be stable under the action of $G$. It is well known that the natural map
\[
G-\wis{Hilb}~\C^d \rOnto X = \C^d/G \]
is a minimal resolution if $d=2$ and if $d=3$ it is often a crepant resolution, for example whenever $G$ is Abelian, see \cite{CrawNotes} for more details. In all cases where $G-\wis{Hilb}~\C^d$ is a desingularization we have again that the corresponding open subvariety $\wis{rep}_{\alpha}^{\theta-semist}~\Lambda$ is smooth. For, in this case the quotient map
\[
\wis{rep}_{\alpha}^{\theta-semist}~\Lambda = \wis{rep}_{\alpha}^{\theta-stable}~\Lambda \rOnto \wis{moduli}^{\theta}_{\alpha}~\Lambda = G-\wis{Hilb}~\C^d \]
is a principal $\wis{PGL}(\alpha)$-fibration and as the base space is smooth by assumption so is the top space.

As we didn't find explicit non-Abelian examples for $\C^3$ in the literature, we include the following simplest example.

Let $A_{4}$ be the alternating group of $12$ elements acting on three dimensional space $\C^3$ via the matrices
\[
A_{4} = \langle~s=\begin{bmatrix}1&0&0\\0&-1&0\\0&0&-1\end{bmatrix}, t=\begin{bmatrix}-1&0&0\\0&-1&0\\0&0&1\end{bmatrix},
r=\begin{bmatrix}0&1&0\\0&0&1\\ 1& 0&0\end{bmatrix}~\rangle \]
the corresponding quotient singularity $\C^3/A_{4}$ has coordinate ring
\[
\C[x,y,z]^{A_{4}}=\C[A(x,y,z),B(x,y,z),C(x,y,z),D(x,y,z)].
\]
with
\[
\begin{cases}
A(x,y,z)=xyz,\\
B(x,y,z)=x^2+y^2+z^2,\\
C(x,y,z)=x^2y^2+y^2z^2+z^2x^2,\\
D(x,y,z)=x^4y^2+y^4z^2+z^4x^2.
\end{cases}
\]
$A$, $B$, $C$ and $D$ obey the relation
\[
D^2+C^3-BCD+A^2(3D-6BC+B^3+9A^2)=0,
\]
whence the quotient singularity $\C^3/A_{4}$ is a hypersurface in $\C^4$.

The character table of the group $A_{4}$ is given by
\[
\begin{array}{c|ccccc}
A_{4} & 1 & \left[\begin{smallmatrix}*&0&0\\0&*&0\\0&0&* \end{smallmatrix}\right] & \left[\begin{smallmatrix}0&*&0\\0&0&*\\ *& 0&0 \end{smallmatrix}\right] & \left[\begin{smallmatrix}0&0&*\\ *&0&0\\ 0&*&0 \end{smallmatrix}\right] \\
& & & & \\
\hline
& & & & \\
V_0 & 1 & 1 & 1 & 1 \\
V_1 & 1 & 1 & \rho & \rho^2\\
V_2 & 1 & 1 & \rho^2 & \rho \\
V_3 & 3 & -1 & 0 & 0
\end{array}
\]
where $\rho$ is a primitive $3$-rd root of unity and therefore
the regular representation is $R=V_0\oplus V_1 \oplus V_2 \oplus V_3^{(1)} \oplus V_3^{(2)} \oplus V_3^{(3)}$. From the character table we deduce the isomorphisms of $A_{4}$-representations
\begin{align*}
&V_3\otimes V_0=V_3\otimes V_1=V_3\otimes V_2=V_3\\
&V_3\otimes V_3=V_0\oplus V_1 \oplus V_2 \oplus V_3 \oplus V_3
\end{align*}
whence the McKay quiver is of the following shape
$$
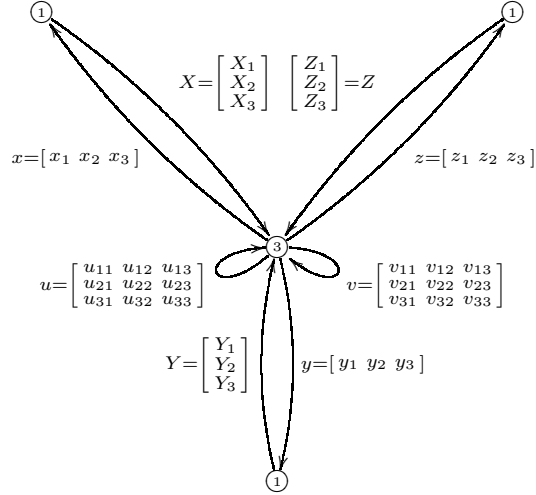
\vcenter{
\xymatrix@=2.8cm{
\vtx{1}\ar@/^/[dr]^{X=\left[\begin{smallmatrix}
X_1\\ X_2\\ X_3
\end{smallmatrix}\right]}  
&
&\vtx{1}\ar@/_/[dl]_{\left[\begin{smallmatrix}
Z_1\\ Z_2\\ Z_3
\end{smallmatrix}\right]=Z}
\\
&\vtx{3}\ar@/^/[ul]^{x=\left[\begin{smallmatrix}
x_1 & x_2 & x_3
\end{smallmatrix}\right]}\ar@/^/[d]^{y=\left[\begin{smallmatrix}
y_1 & y_2 & y_3
\end{smallmatrix}\right]}\ar@/_/[ur]_{z=\left[\begin{smallmatrix}
z_1 & z_2 & z_3
\end{smallmatrix}\right]}\ar@(ld,l)^{u=\left[\begin{smallmatrix}
u_{11} & u_{12} & u_{13}\\
u_{21} & u_{22} & u_{23}\\
u_{31} & u_{32} & u_{33}
\end{smallmatrix}\right]}\ar@(r,rd)^{v=\left[\begin{smallmatrix}
v_{11} & v_{12} & v_{13}\\
v_{21} & v_{22} & v_{23}\\
v_{31} & v_{32} & v_{33}
\end{smallmatrix}\right]}
\\
&\vtx{1}\ar@/^/[u]^{Y=\left[\begin{smallmatrix}
Y_1\\ Y_2\\ Y_3
\end{smallmatrix}\right]} }
}
$$
Denoting
$
V_0=\C v_0, V_1=\C v_1, V_2=\C v_2
$
and
$
V_3^{(i)}=\C e_1^{(i)} + \C e_2^{(i)} + \C e_3^{(i)},
$
we construct a $G$-equivariant basis for
\begin{align*}
V \otimes R &= V_3 \oplus V_3 \oplus V_3 \oplus (V_0\oplus V_1\oplus V_2\oplus V_3\oplus V_3)\\ &\quad \oplus (V_0\oplus V_1\oplus V_2\oplus V_3\oplus V_3) \oplus (V_0\oplus V_1\oplus V_2\oplus V_3\oplus V_3)
\end{align*}
determined by
\begin{alignat*}{2}
V \otimes V_0&=\C(e_1 \otimes v_0) + \C(e_2 \otimes v_0) + \C(e_3 \otimes v_0)\\
V \otimes V_1&=\C(\rho^2 e_1 \otimes v_1) + \C(e_2 \otimes v_1) + \C(\rho e_3 \otimes v_1)\\
V \otimes V_2&=\C(\rho e_1 \otimes v_2) + \C(e_2 \otimes v_2) + \C(\rho^2 e_3 \otimes v_2)\\
V \otimes V_3^{(i)}&=\C(e_1 \otimes e_1^{(i)}) + \C(e_2\otimes e_2^{(i)}) + \C(e_3 \otimes e_3^{(i)})&\qquad &(V_0)\\
& \quad +\C(\rho^2 e_1 \otimes e_1^{(i)}) + \C(\rho e_2 \otimes e_2^{(i)}) + \C( e_3 \otimes e_3^{(i)})&\qquad &(V_1)\\
& \quad +\C(\rho e_1 \otimes e_1^{(i)}) + \C(\rho^2 e_2 \otimes e_2^{(i)}) + \C( e_3 \otimes e_3^{(i)})&\qquad &(V_2)\\
& \quad +\C(e_2 \otimes e_3) + \C(e_3 \otimes e_1) + \C(e_1 \otimes e_2)&\qquad &(V_3 \sim u)\\
& \quad +\C(e_1 \otimes v_0) + \C(e_2 \otimes v_0) + \C(e_3 \otimes v_0)&\qquad &(V_3 \sim v)
\end{alignat*}
With respect to this basis we obtain the following three $12 \times 12$ matrices
\[
P=\left[\begin{smallmatrix}
0 & 0& 0 & x_1 & 0 & 0 & x_2 & 0 & 0& x_3 & 0 &  0\\
0 & 0& 0 & \rho^2y_1 & 0 & 0 & \rho^2y_2 & 0 & 0& \rho^2y_3 & 0 &  0\\
0 & 0& 0 & \rho z_1 & 0 & 0 & \rho z_2 & 0 & 0& \rho z_3 & 0 &  0\\
X_1 & \rho^2Y_1& \rho Z_1 &0 & 0& 0 &0 & 0&0 &0 & 0& 0 \\
0 & 0& 0 &0 & 0& u_{11} &0 & 0& u_{12} &0 & 0& u_{13} \\
0& 0 &0 & 0& v_{11} &0 & 0& v_{12} &0 & 0& v_{13} &0 \\
X_2 & \rho^2Y_2& \rho Z_2 &0 & 0& 0 &0 & 0&0 &0 & 0& 0 \\
0 & 0& 0 &0 & 0& u_{21} &0 & 0& u_{22} &0 & 0& u_{23} \\
0& 0 &0 & 0& v21 &0 & 0& v_{22} &0 & 0& v_{23} &0 \\
X_3 & \rho^2Y_3& \rho Z_3 &0 & 0& 0 &0 & 0&0 &0 & 0& 0 \\
0 & 0& 0 &0 & 0& u31 &0 & 0& u_{32} &0 & 0& u_{33} \\
0& 0 &0 & 0& v_{31} &0 & 0& v_{32} &0 & 0& v_{33} &0
\end{smallmatrix}\right], \]
\[
Q=\left[\begin{smallmatrix}
0 &0 & 0& 0 & x_1 & 0 & 0 & x_2 & 0 & 0& x_3 & 0 \\
0 &0 & 0& 0 & \rho y_1 & 0 & 0 & \rho y_2 & 0 & 0& \rho y_3 & 0\\
0 &0 & 0& 0 & \rho^2 z_1 & 0 & 0 & \rho^2 z_2 & 0 & 0& \rho^2 z_3 & 0 \\
0& 0&0 &0 & 0& v_{11} &0 & 0& v_{12} &0 & 0& v_{13} \\
X_1 & Y_1& Z_1 &0 & 0& 0 &0 & 0&0 &0 & 0& 0 \\
 0 &0 & 0& u_{11} &0 & 0& u_{12} &0 & 0& u_{13} & 0&0\\
0& 0&0 &0 & 0& v_{21} &0 & 0& v_{22} &0 & 0& v_{23} \\
X_2 & Y_2& Z_2 &0 & 0& 0 &0 & 0&0 &0 & 0& 0 \\
 0 &0 & 0& u_{21} &0 & 0& u_{22} &0 & 0& u_{23} &0 &0\\
0& 0 &0&0 & 0& v_{31} &0 & 0& v_{32} &0 & 0& v_{33} \\
X_3 & Y_3& Z_3 &0 & 0& 0 &0 & 0&0 &0 & 0& 0 \\
 0 &0 & 0& u_{31} &0 & 0& u_{32} &0 & 0& u_{33} &0 &0
\end{smallmatrix}\right],
\]
\[
R=\left[\begin{smallmatrix}
0&0&0 & 0& 0 & x_1 & 0 & 0 & x_2 & 0 & 0& x_3\\
0&0&0 & 0& 0 & y_1 & 0 & 0 & y_2 & 0 & 0& y_3\\
0&0&0 & 0& 0 & z_1 & 0 & 0 & z_2 & 0 & 0& z_3 \\
0& 0 &0 & 0& u_{11} &0 & 0& u_{12} &0 & 0& u_{13}  &0\\
0 &0 & 0& v_{11} &0 & 0& v_{12} &0 & 0& v_{13} &0  & 0\\
X_1 & \rho Y_1& \rho^2 Z_1 &0 & 0& 0 &0 & 0&0 &0 & 0& 0 \\
0& 0 &0 & 0& u_{21} &0 & 0& u_{22} &0 & 0& u_{23} & 0 \\
0 &0 & 0& v_{21} &0 & 0& v_{22} &0 & 0& v_{23} &0 & 0\\
X_2 & \rho Y_2& \rho^2 Z_2 &0 & 0& 0 &0 & 0&0 &0 & 0& 0 \\
0& 0 &0 & 0& u_{31} &0 & 0& u_{32} &0 & 0& u_{33} &0\\
0 &0 & 0& v_{31} &0 & 0& v_{32} &0 & 0& v_{33} &0 & 0\\
X_3 & \rho Y_3& \rho^2 Z_3 &0 & 0& 0 &0 & 0&0 &0 & 0& 0
\end{smallmatrix}\right].
\]
Setting the three commutators equal to 0, we obtain the constraints:
\begin{gather*}
x(u-v)=0, \quad y(u-\rho^2 v)=0, \quad z(u-\rho v)=0,\\
(u-v)X=0, \quad (u-\rho^2 v)Y=0, \quad (u-\rho v)Z=0,
u^2=Xx+Yy+Zz,\\
v^2=Xx+\rho^2Yy+\rho Zz.
\end{gather*}
recovering the result obtained in \cite{BG98}.
\end{example}

\section{Partial desingularizations}

In the previous section we have seen that in many cases of current interest one associates to a singularity $\mathfrak{m}$ an $R$-order $\Lambda$ and a stability structure $\theta$ for the dimension vector $\alpha$ such that $\wis{trep}_n~\Lambda = \wis{GL}_n \times^{\wis{GL}(\alpha)} \wis{rep}_{\alpha}~\Lambda$, such that the Zariski open subset
\[
\wis{rep}_{\alpha}^{\theta-semist}~\Lambda \]
of $\theta$-semistable representations is a smooth variety. If this is the case we will call $(\Lambda,\alpha,\theta)$ a {\em good $\mathfrak{m}$-setting}. In this section we will prove that to a good $\mathfrak{m}$-setting one associates a non-commutative desingularization of $\mathfrak{m}$ and a partial commutative desingularization with excellent control on the remaining singularities. We will sketch the procedure in general and then give an explicit description in case $\Lambda$ is a {\em quiver-order}. That is, if 
\[
\Lambda \simeq \int_{\alpha}~\frac{\C Q}{I} \]
for some dimension vector $\alpha$ such that $\wis{rep}_{\alpha}~\C Q/I$ contains (a Zariski open subset of) simple representations and where $\int_{\alpha} \C Q/I$ denotes the algebra of $\wis{GL}_n$-equivariant maps 
\[
\wis{GL}_n \times^{\wis{GL}(\alpha)} \wis{rep}_{\alpha}~\frac{\C Q}{I} \rTo M_n(\C) \]
if $n$ is the total dimension of $\alpha$.

If $(\Lambda,\alpha,\theta)$ is a good $\mathfrak{m}$-setting we have the diagram explained in the previous section
\[
\xymatrix@R=45pt@C=45pt{
\wis{rep}_{\alpha}^{\theta-semist}~\Lambda \ar@{->>}[d]_q \ar@{->>}[rd]^{q_c} \\
\wis{moduli}_{\alpha}^{\theta}~\Lambda \ar@{->>}[r]^{\pi} & X = \wis{spec}~R
}
\]
where $q$ is the algebraic quotient map and $\pi$ is a projective birational map. To $q$ we will assign a sheaf of smooth orders $\mathcal{A}$ on $\wis{moduli}^{\theta}_{\alpha}~\Lambda$. Let $\cup_D~X_D$ be a Zariski open covering by affine normal varieties of the moduli space $\wis{moduli}_{\alpha}^{\theta}~\Lambda$, then each $X_D$ determines a smooth order $\Lambda_D$ defined by taking the algebra of $\wis{GL}_n$-equivariant maps
\[
\wis{GL}_n \times^{\wis{GL}(\alpha)} q^{-1}(X_D) \rTo M_n(\C) \]
for which $q^{-1}(X_D) \simeq \wis{rep}_{\alpha}~\Lambda_D$. Remark that as $q^{-1}(X_D)$ is a smooth $\wis{GL}(\alpha)$-affine variety, we have that
\[
\wis{trep}_n~\Lambda_D = \wis{GL}_n \times^{\wis{GL}(\alpha)} q^{-1}(X_D) \]
is a smooth $\wis{GL}_n$-variety and therefore $\Lambda_D$ is indeed a smooth order. Taking as sections
\[
\Gamma(X_D,\mathcal{A}) = \Lambda_D, 
\]
we obtain a sheaf of smooth orders on $\wis{moduli}_{\alpha}^{\theta}~\Lambda$. We will construct the orders $\Lambda_D$ explicitly if $\Lambda$ is a quiver-order $\int_{\alpha}~\C Q/I$.

Because $\wis{moduli}^{\theta}_{\alpha}~\Lambda = \wis{proj}~\C[\wis{rep}_{\alpha}~\Lambda]^{\wis{GL}(\alpha),\theta}$ we need control on the generators of all $\theta$-semi-invariants. Such a generating set was found by Aidan Schofield and Michel Van den Bergh in \cite{SchofVdB}: {\em determinantal semi-invariants}. In order to define them we have to introduce some notation first.

Reorder the vertices in $Q$ such that the entries of $\theta$ are separated in three strings
\[
\theta = (\underbrace{t_1,\hdots,t_i}_{> 0},\underbrace{t_{i+1},\hdots,t_j}_{=0},\underbrace{t_{j+1},\hdots,t_k}_{< 0}) \]
and let $\theta$ be such that $\theta.\alpha = 0$. Fix a nonzero weight $l \in \N$ and take arbitrary natural numbers $\{ l_{i+1},\hdots,l_j \}$. 

Consider a rectangular matrix $L$ with
\begin{itemize}
\item{$lt_1+\hdots+lt_i+l_{i+1} + \hdots + l_j$ rows and}
\item{$l_{i+1} + \hdots + l_j - l t_{j+1} - \hdots - l t_k$ columns}
\end{itemize}
\[
L = \quad \begin{array}{cc||c|c|c|c|c|c}
& & \overbrace{}^{l_{i+1}} & \hdots & \overbrace{}^{l_j} & \overbrace{}^{-lt_{j+1}} & \hdots & \overbrace{}^{-lt_k} \\
\hline \hline 
lt_1&  \{ & L_{1,i+1} & & L_{1,j} & L_{1,j+1} & & L_{1,k} \\
\hline
& \vdots & & & & & \\
\hline 
lt_i & \{ & L_{i,i+1} & & L_{i,j} & L_{i,j+1} & & L_{i,k} \\
\hline
l_{i+1} & \{ & L_{i+1,i+1} & & L_{i+1,j} & L_{i+1,j+1} & & L_{i+1,k} \\
\hline
& \vdots & & & & & \\
\hline
l_j & \{ & L_{j,i+1} & & L_{j,j}  & L_{j,j+1} & & L_{j,k} 
\end{array}
\]

in which each entry of $L_{r,c}$ is a linear combination of oriented paths in the quiver $Q$ with starting vertex $v_c$ and ending vertex $v_r$.

The relevance of this is that we can evaluate $L$ at any representation $V \in \wis{rep}_{\alpha}~\Lambda$ and obtain a {\em square matrix} $L(V)$ as $\theta.\alpha = 0$. More precisely, if $V_i$ is the vertex-space of $V$ at vertex $v_i$ (that is, $V_i$ has dimension $e_i$), then evaluating $L$ at $V$ gives a linear map
\[
\xymatrix@R=40pt@C=45pt{
V_{i+1}^{\oplus l_{i+1}} \oplus \hdots \oplus V_j^{\oplus l_j} \oplus V_{j+1}^{\oplus -lt_{j+1}} \oplus \hdots \oplus V_k^{\oplus -lt_k} \ar[d]^{L(V)}\\
V_1^{\oplus lt_1} \oplus \hdots \oplus V_i^{\oplus lt_i} \oplus V_{i+1}^{\oplus l_{i+1}} \oplus \hdots \oplus V_j^{\oplus l_j}
}
\]
and $L(V)$ is a square $N \times N$ matrix where
\[
l_{i+1} + \hdots + l_j - lt_{j+1} - \hdots - lt_k = N = lt_1 + \hdots + lt_i + l_{i+1} + \hdots + l_j. \]
So we can consider $D(V) = \wis{det} L(V)$ and verify that $D$ is a $\wis{GL}(\alpha)$-semi-invariant polynomial on $\wis{rep}_{\alpha}~\Lambda$ of weight $\chi_{\theta}^l$. The result of \cite{SchofVdB} asserts that these {\em determinantal semi-invariants} are algebra generators of the graded algebra
\[
\C[\wis{rep}_{\alpha}~\Lambda]^{\wis{GL}(\alpha),\theta}. \]
Observe that this result is to semi-invariants what the result of \cite{LBProcesi} is to invariants. In fact, one can deduce the latter from the first. 

We have seen that a representation $V \in \wis{rep}_{\alpha}~\Lambda$ is $\theta$-semistable if and only if some semi-invariant of weight $\chi_{\theta}^l$ for some $l$ is non-zero on it. This proves

\begin{theorem} The Zariski open subset of $\theta$-semistable $\alpha$-dimensional $\Lambda$-representations can be covered by affine $\wis{GL}(\alpha)$-stable open subsets
\[
\wis{rep}^{\theta-semist}_{\alpha}~\Lambda = \bigcup_D \{ V~|~D(V) = \wis{det} L(V) \not= 0 \} \]
and hence the moduli space can also be covered by affine open subsets
\[
\wis{moduli}^{\theta}_{\alpha}~\Lambda = \bigcup_D~X_D \]
where
$
X_D = \{ [V] \in \wis{moduli}^{\theta}_{\alpha}~\Lambda~|~D(V)=\wis{det} L(V) \not= 0 \}
$.
\end{theorem}

Analogous to the rectangular matrix $L$ we define a rectangular matrix $N$ with
\begin{itemize}
\item{$lt_1+\hdots+lt_i+l_{i+1} + \hdots + l_j$ columns and}
\item{$l_{i+1} + \hdots + l_j - l t_{j+1} - \hdots - l t_k$ rows}
\end{itemize}
\[
N = \quad \begin{array}{cc||c|c|c|c|c|c}
& & \overbrace{}^{l t_1} & \hdots & \overbrace{}^{l t_i} & \overbrace{}^{l_{i+1}} & \hdots & \overbrace{}^{l_j} \\
\hline \hline
l_{i+1} &  \{ & N_{i+1,1} & & N_{i+1,i} & N_{i+1,i+1} & & N_{i+1,j} \\
\hline
& \vdots & & & & & \\
\hline
l_j & \{ & N_{j,1} & & N_{j,i} & N_{j,i+1} & & N_{j,j} \\
\hline
-lt_{j+1}  & \{ & N_{j+1,1} & & N_{j+1,i} & N_{j+1,i+1} & & N_{j+1,j} \\
\hline
& \vdots & & & & & \\
\hline
-l t_k & \{ & N_{k,1} & & N_{k,i}  & N_{k,i+1} & & N_{k,j} 
\end{array}
\]
filled with new variables and define an {\em extended quiver} $Q_D$ where we adjoin for each entry in $N_{r,c}$ an additional arrow from $v_c$ to $v_r$ and denote it with the corresponding variable from $N$.

Let $I_1$ (resp. $I_2$) be the set of relations in $\C Q_D$ determined from the matrix-equations
{\tiny
\[
N.L = \begin{bmatrix} \boxed{(v_{i+1})_{l_{i+1} } } & & & & & 0 \\
& \ddots & & & & \\
& & \boxed{(v_j)_{l_j}} & & & \\
& & & \boxed{(v_{j+1})_{-lt_{j+1}} }& &  \\
& & & & \ddots & \\
0 & & & & & \boxed{(v_k)_{-lt_k}}
\end{bmatrix}
\]}
respectively
{\tiny
\[
L.N = \begin{bmatrix}
\boxed{(v_1)_{lt_1}} & & & & & 0 \\
& \ddots & & & & \\
& & \boxed{(v_i)_{lt_i} }& & & \\
& & & \boxed{(v_{i+1})_{l_{i+1}} }& & \\
& & & & \ddots & \\
0 & & & & & \boxed{(v_j)_{l_j}}
\end{bmatrix}
\]}
where $(v_i)_{n_j}$ is the square $n_j \times n_j$ matrix with $v_i$ on the diagonal and zeroes elsewhere.
Define a new quiver order
\[
\Lambda_D = \int_{\alpha}~\frac{\C Q_D}{(I,I_1,I_2)}
\]
then $\Lambda_D$ is a $\C[X_D]$-order in $\wis{alg@n}$. In fact, the construction of $\Lambda_D$ is nothing but a universal localization in the category $\wis{alg@}\alpha$, which is the subcategory of $\wis{alg@n}$ consisting of all $S = \underbrace{\C \times \hdots \times \C}_k$-algebras with trace map specified by $\alpha$.

That is, take $P_i = v_i \Lambda$ be the projective right ideal associated to vertex $v_i$, then $L$ determines a $\Lambda$-module morphism
\[
P = P_{i+1}^{\oplus l_{i+1}} \oplus \hdots \oplus P_k^{\oplus -lt_k} \rTo^{L} P_1^{\oplus lt_1} \oplus \hdots \oplus P_j^{\oplus l_j} = Q. \]
The algebra map $\Lambda \rTo^{\phi} \Lambda_D$ is universal in $\wis{alg@}\alpha$ with respect to $L \otimes \phi$ being invertible, that is, if $\Lambda \rTo^{\psi} B$ is a morphism in $\wis{alg@}\alpha$ such that $L \otimes \psi$ is an isomorphism of right $B$-modules, then there is a unique map in $\wis{alg@}\alpha$ $\Lambda_D \rTo^u B$ such that $\psi = u \circ \phi$. We claim to have the following situation
\[
\xymatrix@R=40pt@C=45pt{
\wis{rep}^{\theta-semist}~\Lambda \ar@{->>}[d]_q & q^{-1}(X_D)  \simeq \wis{rep}_{\alpha}~\Lambda_D \ar@{_(->}[l] \ar@{->>}[d]\\
\wis{moduli}^{\theta}_{\alpha}~\Lambda & X_D \ar@{_(->}[l]
}
\]
which follows from the next lemma.

\begin{lemma} The following statements are equivalent
\begin{enumerate}
\item{$V \in \wis{rep}_{\alpha}^{\theta-semist}~\Lambda$ lies in $q^{-1}(X_D)$, and}
\item{There is a unique extension $\tilde{V}$ of $V$ such that $\tilde{V} \in \wis{rep}_{\alpha}~\Lambda_D$.}
\end{enumerate}
\end{lemma}

\begin{proof}
$1 \Rightarrow 2$ : Because $L(V)$ is invertible we can take $N(V)$ to be its inverse and decompose it into blocks corresponding to the new arrows in $Q_D$. This then defines the unique extension $\tilde{V} \in \wis{rep}_{\alpha}~Q_D$ of $V$. As $\tilde{V}$ satisfies $I$ (because $V$ does) and $I_1$ and $I_2$ (because $N(V) = L(V)^{-1}$) we have that $\tilde{V} \in \wis{rep}_{\alpha}~\Lambda_D$.

$2 \Rightarrow 1$ :  Restrict $\tilde{V}$ to the arrows of $Q$ to get a $V \in \wis{rep}_{\alpha}~Q$. As $\tilde{V}$ (and hence $V$) satisfies $I$, $V \in \wis{rep}_{\alpha}~\Lambda$. Moreover, $V$ is such that $L(V)$ is invertible (this follows because $\tilde{V}$ satisfies $I_1$ and $I_2$). Hence, $D(V) \not= 0$ and because $D$ is a $\theta$-semi-invariant it follows that $V$ is an $\alpha$-dimensional $\theta$-semistable representation of $\Lambda$. An alternative method to see this is as follows. Assume that $V$ is {\em not} $\theta$-semistable and let $V' \subset V$ be a subrepresentation such that $\theta.\wis{dim} V' < 0$. Consider the restriction of the linear map $L(V)$ to the subrepresentation $V'$ and look at the commuting diagram
\[
\xymatrix@R=40pt@C=45pt{
V_{i+1}^{'\oplus l_{i+1}} \oplus \hdots \oplus V_k^{'\oplus -lt_k} \ar[r]^{L(V)|V'} \ar@{^(->}[d] & V_1^{'\oplus lt_1} \oplus \hdots \oplus V_j^{'\oplus l_j} \ar@{^(->}[d] \\
V_{i+1}^{\oplus l_{i+1}} \oplus \hdots \oplus V_k^{\oplus -lt_k} \ar[r]^{L(V)} & V_1^{\oplus lt_1} \oplus \hdots \oplus V_j^{\oplus l_j}
}
\]
As $\theta. \wis{dim} V' < 0$ the top-map must have a kernel which is clearly absurd as we know that $L(V)$ is invertible.
\end{proof}

The universal property of the universal localizations $\Lambda_D$ allows us to glue these orders together into a coherent sheaf on $\wis{moduli}_{\alpha}^{\theta}~\Lambda$. Let $\Lambda_{D_1}$ (resp. $\Lambda_{D_2}$) be the order constructed from a rectangular matrix $L_1$ (resp. $L_2$), then we can construct the direct sum map $L = L_1 \oplus L_2$ for which the corresponding semi-invariant $D=D_1D_2$. As $\Lambda \rTo \Lambda_D$ makes the projective module morphisms associated to $L_1$ and $L_2$ into an isomorphism we have uniquely determined maps in $\wis{alg@}\alpha$
\[
\xymatrix{
& \Lambda_D \\
\Lambda_{D_1} \ar[ur]^{i_1} & & \Lambda_{D_2} \ar[ul]_{i_2}
}
\qquad \text{whence}
\qquad
\xymatrix@C=10pt{
& \wis{rep}_{\alpha}~\Lambda_D \ar[dl]_{i_1^*}  \ar[dr]^{i_2^*}\\
\wis{rep}_{\alpha}~\Lambda_{D_1} & & \wis{rep}_{\alpha}~\Lambda_{D_2} 
}
\]
Because $\wis{rep}_{\alpha}~\Lambda_D = q^{-1}(X_D)$ (and similarly for $D_i$) we have that $i_j^*$ are embeddings as are the $i_j$. This way we can glue the sections $\Gamma(X_{D_1},\Ascr) = \Lambda_{D_1}$ with $\Gamma(X_{D_2},\Ascr) = \Lambda_{D_2}$ over their intersection $X_D = X_{D_1} \cap X_{D_2}$ via the inclusions $i_j$. Hence we get a coherent sheaf of non-commutative algebras $\Ascr$ over $\wis{moduli}^{\theta}_{\alpha}~\Lambda$. Further, by localizing the orders $\Lambda_{D_j}$ at the central element $D$ we have that the algebra morphisms $i_j$ are central extensions, that is satisfying
\[
\Lambda_D = \Lambda_{D_j} Z(\Lambda_D) \]
which implies that we have morphisms between the non-commutative structure sheaves
\[
(\wis{spec}~\Lambda_{D_j},\Oscr^{nc}_{\Lambda_{D_j}}) \rTo (\wis{spec}~\Lambda_D,\Oscr^{nc}_{\Lambda_D}) \]
which allow us to define a non-commutative variety $\wis{spec}~\Ascr$ by gluing the non-commutative structure sheaves of the various $\Lambda_{D_j}$ together. Observe that the central scheme of this non-commutative variety is $\wis{moduli}_{\alpha}^{\theta}~\Lambda$ with its structure sheaf. This concludes the proof of the following result.

\begin{theorem} Let $(\Lambda,\alpha,\theta)$ be a good $\mathfrak{m}$-setting. Then, there is a sheaf of smooth orders $\Ascr$ over the moduli space $\wis{moduli}_{\alpha}^{\theta}~\Lambda$ such that the diagram below is commutative
\[
\xymatrix@R=40pt@C=45pt{
\wis{spec}~\Ascr \ar[d]_c \ar[rd]^{\phi} \\
\wis{moduli}^{\theta}_{\alpha}~\Lambda \ar@{->>}[r]^{\pi} & X = \wis{spec}~R
}
\]
Here, $\wis{spec}~\Ascr$ is a non-commutative variety obtained by gluing affine non-commutative structure sheaves $(\wis{spec}~\Lambda_D,\Oscr^{nc}_{\Lambda_D})$ together and where $c$ is the map which intersects locally a prime ideal of $\Lambda_D$ with its center. Because $\Ascr$ is a sheaf of smooth orders in $\wis{alg@n}$, $\phi$ can be viewed as a {\em non-commutative desingularization} of $X$.

Moreover, if $\theta$ is such that all $\theta$-semistable $\alpha$-dimensional $\Lambda$-representations are actually $\theta$-stable, then $\Ascr$ is a sheaf of Azumaya algebras over $\wis{moduli}^{\theta}_{\alpha}~\Lambda$ and in this case $\pi$ is a commutative desingularization of $X$. If, in addition, also $\alpha$ is an indivisible dimension vector (that is, $gcd(\alpha) = 1$) then
$\Ascr \simeq End~\mathcal{P}$ for some vectorbundle $\mathcal{P}$ of rank $n$ over $\wis{moduli}^{\theta}_{\alpha}~\Lambda$.
\end{theorem}

In general, there may remain singularities in $\wis{moduli}^{\theta}_{\alpha}~\Lambda$ but then have been fully classified in dimensions $\leq 6$ and reduction steps exists which prove that in each dimension there is a finite list of such possible remaining singularities. We will recall these steps briefly, the starting point being the local marked quiver setting $(Q^{\dagger},\alpha)$ associated to a point $\mathfrak{n} \in \wis{moduli}^{\theta}_{\alpha}~\Lambda$. Remark that $\mathfrak{n} \in X_D$ for some $D$ and as $\Lambda_D$ is a smooth order in $\wis{alg@n}$ the defect $\wis{def}_{\mathfrak{n}}~\Lambda_D = 0$ so the local marked quiver setting determines the \'etale local structure of $\wis{moduli}^{\theta}_{\alpha}~\Lambda$ near $\mathfrak{n}$.

The reduction steps below were discovered by R. Bocklandt in his Ph.D. thesis \cite{BocklandtThesis} (see also \cite{Bocklandtpaper}) in which he classifies quiver settings having a regular ring of invariants. These steps were slightly extended in \cite{RBLBVdW} in order to classify central singularities of smooth orders. All reductions are made locally around a vertex in the marked quiver. There are three types of allowed moves 

\par \vskip 3mm
\noindent
{\bf 1.Vertex removal}
Assume we have a marked quiver setting $(Q^{\dagger},\alpha)$ and a vertex $v$ such that the local structure of $(Q^{\dagger},\alpha)$ near $v$ is indicated by the picture on the left below, that is, inside the vertices we have written the components of the dimension vector and the subscripts of an arrow indicate how many such arrows there are in $Q^{\dagger}$ between the indicated vertices.
Define the new marked quiver setting $(Q^{\dagger}_R,\alpha_R)$ obtained by the operation $R^v_V$ which removes the vertex $v$ and composes all arrows through $v$, the dimensions of the other vertices are unchanged :
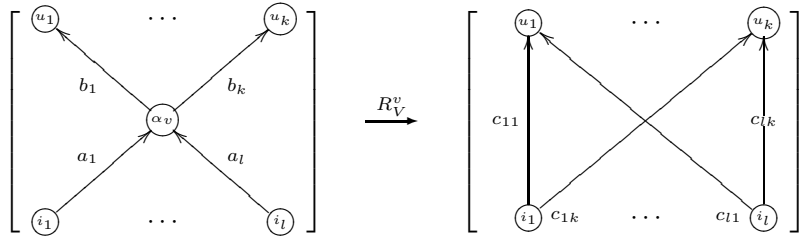
\[
\left[ ~\vcenter{
\xymatrix@=1cm{
\vtx{u_1}&\cdots &\vtx{u_k}\\
&\vtx{\alpha_v}\ar[ul]^{b_1}\ar[ur]_{b_k}&\\
\vtx{i_1}\ar[ur]^{a_1}&\cdots &\vtx{i_l}\ar[ul]_{a_l}}}
~\right] \quad
\rTo^{R^v_V} \quad
\left[~\vcenter{
\xymatrix@=1cm{
\vtx{u_1}&\cdots &\vtx{u_k}\\
&&\\
\vtx{i_1}\ar[uu]^{c_{11}}\ar[uurr]_<<{c_{1k}}&\cdots &\vtx{i_l}\ar[uu]|{c_{lk}}\ar[uull]^<<{c_{l1}}}}
~\right].
\]
where $c_{ij} = a_ib_j$ (observe that some of the incoming and outgoing vertices may be the
same so that one obtains loops in the corresponding vertex). One is allowed to make this reduction step provided either of the following conditions is met
\[
 \chi_Q(\alpha,\epsilon_v) \geq 0 \quad \Leftrightarrow \quad \alpha_v \geq \sum_{j=1}^l a_j i_j \] 
\[
 \chi_Q(\epsilon_v,\alpha) \geq 0\quad \Leftrightarrow \quad \alpha_v \geq \sum_{j=1}^k b_j u_j \]
(observe that if we started off from a marked quiver setting $(Q^{\dagger},\alpha)$ coming from an order, then these inequalities must actually be equalities).

\par \vskip 3mm
\noindent
{\bf 2. loop removal}
If $v$ is a vertex with vertex-dimension $\alpha_v = 1$ and having $k \geq 1$ loops, then let $(Q^{\dagger}_R,\alpha_R)$ be the marked quiver setting obtained by the loop removal operation $R^v_l$
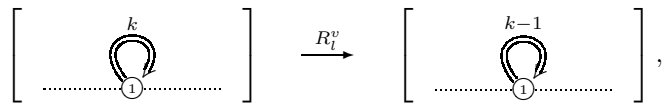
\[
\left[~\vcenter{
\xymatrix@=1cm{
&\vtx{1}\ar@{..}[r]\ar@{..}[l]\ar@(lu,ru)@{=>}^k&}}
~\right]\quad \rTo^{R^v_l} \quad
\left[~\vcenter{
\xymatrix@=1cm{
&\vtx{1}\ar@{..}[r]\ar@{..}[l]\ar@(lu,ru)@{=>}^{k-1}&}}
~\right],\]
removing one loop in $v$ and keeping the same dimension vector. 

\par \vskip 3mm
\noindent
{\bf 3. Loop removal}
If the local situation in $v$ is such that there is exactly one (marked) loop in $v$, the dimension vector in $v$ is $k \geq 2$ and there is exactly one arrow leaving $v$ and this to a vertex with dimension vector $1$, then one is allowed to make the reduction $R^v_L$ indicated below
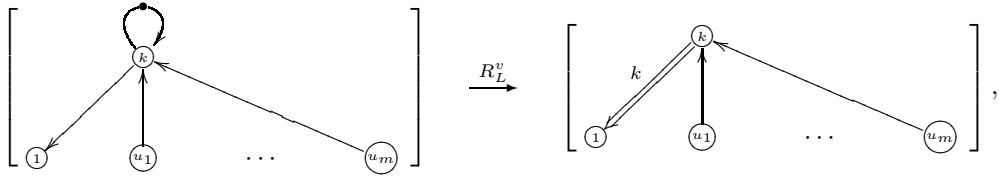
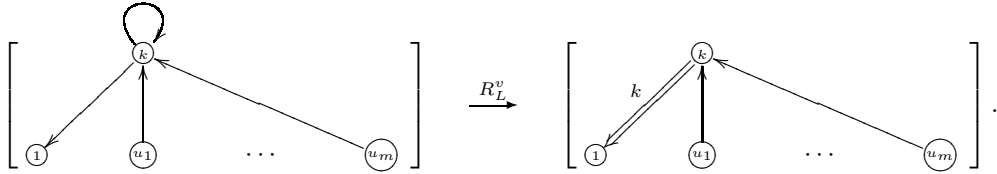
\[
\left[~\vcenter{
\xymatrix@=1cm{
&\vtx{k}\ar[dl]\ar@(lu,ru)|{\bullet}&&\\
\vtx{1}&\vtx{u_1}\ar[u]&\cdots &\vtx{u_m}\ar[ull]}}
~\right]\quad \rTo^{R^v_L} \quad
\left[~\vcenter{
\xymatrix@=1cm{
&\vtx{k}\ar@2[dl]_{k}&&\\
\vtx{1}&\vtx{u_1}\ar[u]&\cdots &\vtx{u_m}\ar[ull]}}
~\right],
\]
\vspace{.5cm}
\[
\left[~\vcenter{
\xymatrix@=1cm{
&\vtx{k}\ar[dl]\ar@(lu,ru)&&\\
\vtx{1}&\vtx{u_1}\ar[u]&\cdots &\vtx{u_m}\ar[ull]}}
~\right]\quad \rTo^{R^v_L} \quad
\left[~\vcenter{
\xymatrix@=1cm{
&\vtx{k}\ar@2[dl]_k&&\\
\vtx{1}&\vtx{u_1}\ar[u]&\cdots &\vtx{u_m}\ar[ull]}}
~\right].
\]

Similarly, if there is one (marked) loop in $v$ and $\alpha_v = k \geq 2$ and there is only one arrow arriving at $v$ coming from a vertex of dimension vector $1$, then one is allowed to make the reduction $R^v_L$
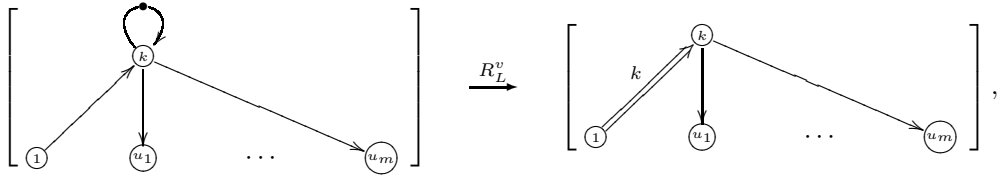
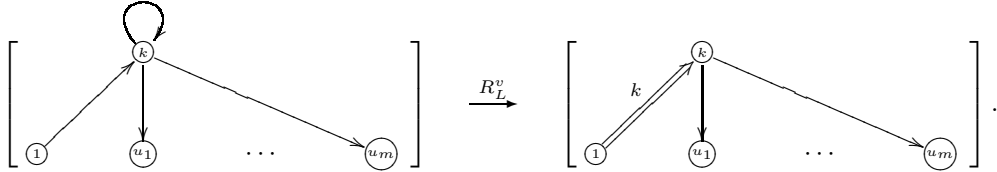
\[
\left[~\vcenter{
\xymatrix@=1cm{
&\vtx{k}\ar[d]\ar[drr]\ar@(lu,ru)|{\bullet}&&\\
\vtx{1}\ar[ur]&\vtx{u_1}&\cdots &\vtx{u_m}}}
~\right]\quad \rTo^{R^v_L} \quad
\left[~\vcenter{
\xymatrix@=1cm{
&\vtx{k}\ar[d]\ar[drr]&&\\
\vtx{1}\ar@2[ur]^k&\vtx{u_1}&\cdots &\vtx{u_m}}}
~\right], \]
\vspace{.5cm}
\[
\left[~\vcenter{
\xymatrix@=1cm{
&\vtx{k}\ar[d]\ar[drr]\ar@(lu,ru)&&\\
\vtx{1}\ar[ur]&\vtx{u_1}&\cdots &\vtx{u_m}}}
~\right]\quad \rTo^{R^v_L} \quad
\left[~\vcenter{
\xymatrix@=1cm{
&\vtx{k}\ar[d]\ar[drr]&&\\
\vtx{1}\ar@2[ur]^k&\vtx{u_1}&\cdots &\vtx{u_m}}}
~\right].
\]
The relevance of these reduction rules is that if
\[
(Q^{\dagger}_1,\alpha_1) \rightsquigarrow (Q^{\dagger}_2,\alpha_2) \]
is a sequence of legal reductions, then 
\[
\C[\wis{rep}_{\alpha_1}~Q^{\dagger}_1]^{\wis{GL}(\alpha_1)} \simeq \C[\wis{rep}_{\alpha_2}~Q^{\dagger}_2]^{\wis{GL}(\alpha_2)}[y_1,\hdots,y_z] \]
where $z$ is the sum of all loops removed in $R^v_l$ reductions plus the sum of $\alpha_v$ for each reduction step $R^v_L$ involving a genuine loop and the sum of $\alpha_v - 1$ for each reduction step $R^v_L$ involving a marked loop. That is, marked quiver settings which belong to the same reduction tree have smooth equivalent invariant rings.

\begin{theorem} Let $(Q^{\dagger},\alpha)$ be a marked quiver setting, then there is a unique reduced setting (that is, having no further admissible reduction steps)
$(Q^{\dagger}_0,\alpha_0)$ for which there exists a reduction procedure
\[
(Q^{\dagger},\alpha) \rightsquigarrow (Q^{\dagger}_0,\alpha_0).. \]
We will denote this unique setting by $Z(Q^{\dagger},\alpha)$.
\end{theorem}

The following result is a slight adaptation of Bocklandt's main result \cite{Bocklandtpaper}.

\begin{theorem} Let $(Q^{\dagger}_{\mathfrak{n}},\alpha_{\mathfrak{n}})$ be the local marked quiver setting of $\mathfrak{n} \in \wis{moduli}^{\theta}_{\alpha}~\Lambda$. Then, $\mathfrak{n}$ is a smooth point  if and only if the unique associated reduced setting
\[
Z(Q^{\dagger}_{\mathfrak{n}},\alpha_{\mathfrak{n}}) \in \{~\xymatrix{\vtx{k}} \qquad \xymatrix{\vtx{k} \ar@(ul,ur)}  \qquad  \xymatrix{\vtx{k} \ar@(ul,ur)|{\bullet}}  \quad~\qquad  \xymatrix{\vtx{2} \ar@(dl,ul) \ar@(dr,ur)} \qquad~\quad~\quad~\qquad~\qquad~\qquad~\qquad~\qquad\xymatrix{\vtx{2} \ar@(dl,ul) \ar@(dr,ur)|{\bullet}}~\quad \xymatrix{\vtx{2} \ar@(dl,ul)|{\bullet} \ar@(dr,ur)|{\bullet}} \qquad~\}.
\]
The Azumaya points are such that $Z(Q^{\dagger}_{\mathfrak{n}},\alpha_{\mathfrak{n}}) = \xymatrix{\vtx{1}}$ hence the singular locus of $\wis{moduli}^{\theta}_{\alpha}~\Lambda$ is contained in the ramification locus of $\mathcal{A}$ but may be strictly smaller.
\end{theorem}

To classify the central singularities of smooth orders we may reduce to zero-settings $(Q^{\dagger},\alpha) = Z(Q^{\dagger},\alpha)$. For such a setting we have for all vertices $v_i$ the inequalities
\[
\chi_Q(\alpha,\delta_i) < 0 \qquad \text{and} \qquad \chi_Q(\delta_i,\alpha) < 0 \]
and the dimension of the central variety can be computed from the Euler-form $\chi_Q$. This gives us an estimate of $d = \wis{dim}~X = \wis{dim}~\wis{moduli}^{\theta}_{\alpha}~\Lambda$ which is very efficient to classify the singularities in low dimensions.

\begin{theorem} \label{counting} Let $(Q^{\dagger},\alpha) = Z(Q^{\dagger},\alpha)$ be a reduced setting on $k \geq 2$ vertices. Then,
\[
\wis{dim}~X  \geq 1 + \sum_{\xymatrix@=1cm{ \vtx{a} }}^{a \geq 1} a + 
\sum_{\xymatrix@=1cm{ \vtx{a}\ar@(ul,dl)|{\bullet} }}^{a > 1}(2a-1) +
 \sum_{\xymatrix@=1cm{ \vtx{a}\ar@(ul,dl)}}^{a > 1}(2a) + \sum_{\xymatrix@=1cm{ \vtx{a}\ar@(ul,dl)|{\bullet}\ar@(ur,dr)|{\bullet}}}^{a > 1} (a^2+a-2) + \]
\[
\sum_{\xymatrix@=1cm{ \vtx{a}\ar@(ul,dl)|{\bullet}\ar@(ur,dr)}}^{a > 1} (a^2+a-1) +
\sum_{\xymatrix@=1cm{ \vtx{a}\ar@(ul,dl)\ar@(ur,dr)}}^{a > 1} (a^2+a) + \hdots +
\sum_{\xymatrix@=1cm{ \vtx{a}\ar@(ul,dl)|{\bullet}_{k}\ar@(ur,dr)^{l}}}^{a > 1} ((k+l-1)a^2+a-k) + \hdots
\]
In this sum the contribution of a vertex $v$ with $\alpha_v = a$ is determined by the number of
(marked) loops in $v$. By the reduction steps (marked) loops only occur at vertices where
$\alpha_v > 1$.
\end{theorem}

For example, this shows that there are no central singularities in dimension $d=2$ and that for $d=3$ the only reduced singular setting is 
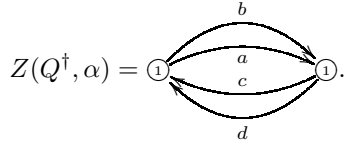
\[
Z(Q^{\dagger},\alpha) = \xymatrix{\vtx{1} \ar@/^2ex/[rr]_a \ar@/^4ex/[rr]^b & & \vtx{1} \ar@/^2ex/[ll]_c \ar@/^4ex/[ll]^d}. \]
The ring of polynomial invariants $R^{\alpha}_{Q^{\dagger}}$ is generated by traces along oriented cycles in $Q^{\dagger}$ so is generated by the invariants
\[
x = ac, \quad y = ad, \quad u = bc \quad \text{and} \quad v = bd~\qquad \text{whence}~\qquad~R^{\alpha}_{Q^{\dagger}} \simeq \frac{\C[x,y,u,v]}{(xy-uv)}. \]
Hence, the only \'etale type of central singularity in dimension three is the {\em conifold singularity}.

\begin{example}[dimension $d=4$]
If $(Q^{\dagger},\alpha)$ is a reduced setting for dimension $4$ then $Q^{\dagger}$ can have at most three vertices. If there is just one, its dimension must be $1$ (smooth setting) or $2$ in which case the only new type is
\[
Z(Q^{\dagger},\alpha) = \qquad \xymatrix{\vtx{2} \ar@(ul,dl) \ar@(ur,dr)|{\bullet}} \]
which is again a smooth setting.
If there are two vertices, both must have dimension $1$ and have at least two incoming and outgoing arrows as in the previous example. The only new type that occurs is
\[
Z(Q^{\dagger},\alpha) = \xymatrix{ \vtx{1} \ar@/^/[rr] \ar@/^3ex/[rr] & & \vtx{1} \ar@/^/[ll] \ar@/^2ex/[ll] \ar@/^3ex/[ll]} \]

\par \vskip 2mm
\noindent
for which one calculates as before the ring of invariants to be
\[
R^{\alpha}_{Q^{\dagger}} = \frac{\C[a,b,c,d,e,f]}{(ae-bd,af-cd,bf-ce)}. \]
If there are three vertices all must have dimension $1$ and each vertex must have at least two incoming and two outgoing arrows. There are just two such possibilities in dimension $4$
\[
Z(Q^{\dagger},\alpha) \in \left\{~\vcenter{\xymatrix{\vtx{1}\ar@/^/[rr]\ar@/^/[rd]&&\vtx{1}\ar@/^/[ll]\ar@/^/[ld]\\
&\vtx{1}\ar@/^/[ru]\ar@/^/[lu]&}}  \qquad \vcenter{ 
\xymatrix{\vtx{1}\ar@2@/^/[rr]&&\vtx{1}\ar@2@/^/[ld]\\
&\vtx{1}\ar@2@/^/[lu]}}~ \right\}.
\]
The corresponding rings of polynomial invariants are
\[
R^{\alpha}_{Q^{\dagger}} = \frac{\C[x_1,x_2,x_3,x_4,x_5]}{(x_4x_5-x_1x_2x_3)} \qquad \text{resp.} \qquad
R^{\alpha}_{Q^{\dagger}} = \frac{\C[x_1,x_2,x_3,x_4,y_1,y_2,y_3,y_4]}{R_2} \]
where $R_2$ is the ideal generated by all $2 \times 2$ minors of the matrix
\[
\begin{bmatrix} 
x_1 & x_2 & x_3 & x_4 \\
y_1 & y_2 & y_3 & y_4 \end{bmatrix}
\]
\end{example}

In \cite{RBLBVdW} it was proved that there are exactly ten types of smooth order central singularities in dimension $d=5$ and $53$ in dimension $d=6$.

\section{The conifold algebra}

Quiver-diagrams play an important role in stringtheory as they encode intersection information of so called {\em wrapped $D$-branes} (higher dimensional strings) in  Calabi-Yau manifolds. One of the earliest models, studied by I. R. Klebanov and E. Witten \cite{KlebanovWitten}, was based on the conifold singularity (see previous section). A {\em $D3$-brane} is a three-dimensional (over the real numbers $\R$) submanifold of a Calabi-Yau manifold and as this is a six-dimensional (again over the real numbers) manifold it follows that two $D3$-branes in sufficiently general position intersect each other in a finite number of points. If one wraps two sufficiently general $D3$-branes around a conifold singularity, their intersection data will be encoded in the quiver-diagram
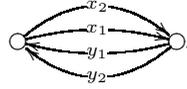
\[
\xymatrix{\vtx{} \ar@/^1ex/[rr]|{x_1} \ar@/^3ex/[rr]|{x_2} & & \vtx{} \ar@/^1ex/[ll]|{y_1} \ar@/^3ex/[ll]|{y_2}}.
\]
Without going into details (for more information see \cite{Berenstein}) one can associate to such a quiver-diagram a non-commutative algebra describing the vacua with respect to a certain {\em super-potential} which is a suitable linear combination of oriented cycles in the quiver-diagram. In the case of two $D3$-branes wrapped around a conifold singularity one obtains :

\begin{definition} The {\em conifold algebra} $\Lambda_c$ is the non-commutative affine $\C$-algebra generated by three non-commuting variables $X,Y$ and $Z$ and satisfying the following relations
\[
\begin{cases}
XZ &= - ZX \\
YZ &= - ZY \\
X^2Y &= YX^2 \\
Y^2X &= XY^2 \\
Z^2 &= 1
\end{cases}
\]
That is, $\Lambda$ has a presentation
\[
\Lambda_c = \frac{\C \langle X,Y,Z \rangle}{(Z^2-1,XZ+ZX,YZ+ZY,[X^2,Y],[Y^2,X])} \]
where $[A,B]=AB-BA$ denotes the commutator.
One sometimes encounters another presentation of $\Lambda_c$ as
\[
\frac{\C \langle X,Y,Z \rangle}{(Z^2-1,XZ+ZX,YZ+ZY,[Z[X,Y],X],[Z[X,Y],Y])} \]
but as $Z$ is a unit, it is easily seen that both presentations give isomorphic $\C$-algebras.
\end{definition}

\begin{proposition} In the conifold algebra $\Lambda_c$ the elements
\[
x = X^2, \qquad y = Y^2 \qquad \text{and} \qquad z = \frac{1}{2}(XY+YX) \]
are algebraically independent central elements and $\Lambda_c$ is a free module over the central subalgebra $C = \C[x,y,z]$ with basis
\[
\Lambda_c = C.1 \oplus C.X \oplus C.Y \oplus C.Z \oplus C.XY \oplus C.XZ \oplus C.YZ \oplus C.XYZ \]
In fact, the conifold algebra is a skew group algebra
\[
\Lambda_c \simeq \C[z,X][Y,\sigma,\delta] \# \Z/2\Z \]
for some automorphism $\sigma$ and $\sigma$-derivation $\delta$. In particular, $\Lambda_c$ is a regular algebra of dimension three.
\end{proposition}

\begin{proof} Consider the subalgebra $S$ of $\Lambda_c$ generated by $X$ and $Y$, that is
\[
S = \frac{\C \langle X,Y \rangle}{([X^2,Y],[Y^2,X])} \]
Then clearly $x$ and $y$ are central elements of $S$ as is $z = \frac{1}{2}(XY+YX)$ because
\[
(XY+YX)X = XYX+YX^2=YXY+X^2Y=X(YX+XY) \]
Now, consider the  \"Ore extension
\[
S' = \C[z,X][Y,\sigma,\delta] \quad \text{with} \quad \sigma(z)=z,\sigma(X)=-X\quad \text{and} \quad \delta(z)=0, \delta(X)=2z \]
This means that $z$ is a central element of $S'$ and that $YX=\sigma(X)Y+\delta(X)=-XY+2z$ whence
the map
\[
S \rTo S' \qquad \text{defined by} \qquad X \mapsto X \quad \text{and} \quad Y \mapsto Y \]
is an isomorphism. By standard results, the {\em center} of $S'$ is equal to
\[
Z(S') = \C[x,y,z] \]
whence the three elements are algebraically independent. Consider the automorphism defined by
$\phi(X) = -X$ and $\phi(Y)=-Y$ on $S$, then the conifold algebra can be written as the {\em skew group ring} 
\[
\Lambda_c \simeq S \# \Z/2\Z \]
As $Z(S) = \C[x,y,z]$ is fixed under $\phi$ the elements $x = x \# 1$, $y = y \# 1$ and $z = z \# 1$ are central in $\Lambda_c$ and as $S'$ is free over $Z(S')$ with basis
\[
S' = Z(S').1 \oplus Z(S').X \oplus Z(S').Y \oplus Z(S').XY \]
the result on freeness of $\Lambda_c$ over $\C[x,y,z]$ follows.
\end{proof}

If $C$ is a commutative $\C$-algebra and if $M_q$ is a {\em symmetric} $m \times m$ matrix with entries in $C$, then we have a {\em bilinear form} on the free $C$-module $V = C \oplus \hdots \oplus C$ of rank $m$ defined by
\[
B_q(v,w) = \begin{bmatrix} v_1 & v_2 &  \hdots & v_m \end{bmatrix}.\begin{bmatrix} b_{11} & b_{12} & \hdots & b_{1n} \\
b_{12} & b_{22} & \hdots & b_{2n} \\
\vdots & \vdots & & \vdots \\
b_{1n} & b_{2n} & \hdots & b_{nn} \end{bmatrix}.\begin{bmatrix} w_1 \\ w_2 \\ \vdots \\ w_m \end{bmatrix}.
\]
The associated {\em Clifford algebra} $Cl_q(V)$ is then the quotient of the {\em tensor algebra} $T_C(V) = C \langle v_1,\hdots,v_m \rangle$ where $\{ v_1,\hdots,v_m \}$ is a basis of the free $C$-module $V$ and the defining relations are
\[
Cl_q(V) = \frac{T_C(V)}{(v \otimes w + w \otimes v - 2B_q(v,w)~:~v,w \in V)} .\]
As an example, the algebra $S \simeq S'$ constructed in the above proof is the Clifford algebra of the binary quadratic form over $C = \C[x,y,z]$
\[
B_q = \begin{bmatrix} x & z \\ z & y \end{bmatrix} \qquad \text{on} \qquad V = C.X \oplus C.Y \]
as $B_q(X,X)=x, B_q(Y,Y)=y$ and $B_q(X,Y) = z$. As the entries of the symmetric variable are independent variables, we call this algebra the {\em generic binary Clifford algebra}, see \cite{LB2x2} for more details and the structure of higher generic Clifford algebras.

\begin{lemma} The conifold algebra $\Lambda_c$ is the {\em Clifford algebra} of a non-degenerate ternary quadratic form over $\C[x,y,z]$.
\end{lemma}

\begin{proof}
Consider the free $C=\C[x,y,z]$-module of rank three $V = C.X \oplus C.Y \oplus C.Z$ and the symmetric $3 \times 3$ matrix
\[
B_q = \begin{bmatrix} x & z & 0 \\ z & y & 0 \\ 0 & 0 & 1 \end{bmatrix}
\]
then it follows that $\Lambda_c \simeq Cl_q(V)$ as $B_q(X,Z)=0, B_q(Y,Z)=0$, $B_q(Z,Z)=0$ and the remaining inproducts are those of $S \simeq S'$ above.
\end{proof}

Whereas $C=\C[x,y,z]$ is a central subalgebra of $\Lambda_c$, the center itself is strictly larger. Take $D=XYZ-YXZ$ and verify that
\begin{eqnarray*}
(XYZ-YXZ)X =& -X(2z-XY)Z+xYZ \\
=& -2zXZ+2xYZ \\
=& xYZ - (2zXZ-YX^2Z) \\
=& X(XYZ-YXZ)
\end{eqnarray*}
and a similar calculation shows that $DY=YD$ and $DZ=ZD$. Moreover, $D \notin \C[x,y,z]$. Indeed, in the description $\Lambda_c \simeq S \# \Z/2\Z$ we have that
\[
\C[x,y,z] \subset S \# 1 \qquad \text{whereas} \qquad D = XYZ-YXZ = (XY-YX) \# Z \in S \# Z. \]
Moreover, we have that $D^2 \in \C[x,y,z]$ because
\[
D^2 = (XYZ-YXZ)^2 = 2z(XY+YX) - 4xy = 4(z^2-xy) \in \C[x,y,z]. \]

\begin{lemma} The center $R_c$ of the conifold algebra $\Lambda_c$ is isomorphic to the coordinate ring of the conifold singularity
\[
R_c \simeq \frac{\C[a,b,c,d]}{(ab-cd)}. \]
\end{lemma}

\begin{proof}
Let $Z$ be the central subalgebra generated by $x,y,z$ and $D$, then a representation of $Z$ is
\[
Z = \frac{\C[x,y,z,D]}{(D^2-4(z^2-xy))} \simeq \frac{\C[a,b,c,d]}{(ab-cd)} \]
where the second isomorphism comes from the following change of coordinates
\[
a = D + 2z,\quad b = D-2z,\quad c=2x \quad \text{and} \quad d = 2y .\]
As a consequence $Z$ is the coordinate ring of the conifold singularity and is in particular integrally closed. As $\Lambda_c$ is a finite module over $Z$ it follows that if $Z \not= R_c$ then the field of fractions $L$ of $R_c$ would be a proper extension of the field of fractions $K$ of $Z$. This can be contradicted using classical results on Clifford algebras over fields. To begin, note that as the ternary form
\[
B_q = \begin{bmatrix} x & z & 0 \\ z & y & 0 \\ 0 & 0 & 1 \end{bmatrix} \]
has square-free determinant $xy-z^2 \notin \C(x,y,z)^{*2}$, the Clifford algebra over the rational field $\C(x,y,z)$ 
\[
\Lambda_c \otimes_{\C[x,y,z]} \C(x,y,z) \]
is a central simple algebra of dimension $4$ over its center $K'$ which is a quadratic field extension of $\C(x,y,z)$ determined by adjoining the square root of the determinant. As $[K : \C(x,y,z)] = 2$ it follows that $K=K'$ and hence also that $K=L$ whence $Z=R_c$.
\end{proof}

Let us relate the non-commutative affine variety $\wis{spec}~\Lambda_c$ with that of the central subalgebra $\wis{spec}~\C[x,y,z] = \Af^3$.

\begin{lemma} Intersecting twosided prime ideals  of $\Lambda_c$ with the central subalgebra $\C[x,y,z]$ determines a continuous map
\[
\wis{spec}~\Lambda_c \rTo^{\phi} \Af^3 \]
with the following fiber information :
\begin{enumerate}
\item{If $\mathfrak{n} \notin \V(xy-z^2)$, then $\phi^{-1}(\mathfrak{n})$ consists of two points.}
\item{If $(x,y,z) \not= \mathfrak{n} \in \V(xy-z^2)$, then $\phi^{-1}(\mathfrak{n})$ consists of one point.}
\item{If $(x,y,z) = \mathfrak{n}$, then $\phi^{-1}(\mathfrak{n})$ consists of two points.}
\end{enumerate}
\end{lemma}

\begin{proof} For $P=(a,b,c) \in \Af^3$  the quotient of $\Lambda_c$ by the extended two-sided ideal $\Lambda_c \mathfrak{n}_P$ is the Clifford algebra $Cl_P$ over $\C$ of the ternary quadratic form
\[
B_P = \begin{bmatrix} a & c & 0 \\ c & b & 0 \\ 0 & 0 & 1 \end{bmatrix} \]
and the elements of $\phi^{-1}(\mathfrak{n}_P)$ are the two-sided maximal ideals of $Cl_P$. We can diagonalize the symmetric matrix, that is there is a base-change matrix $M \in \wis{GL}_3$ such that
\[
M^{\tau}.\begin{bmatrix} a & c & 0 \\ c & b & 0 \\ 0 & 0 & 1 \end{bmatrix}.M = \begin{bmatrix} u & 0 & 0 \\ 0 & v & 0 \\ 0 & 0 & 1 \end{bmatrix} = B_Q \]
(with $uv = ab-c^2$) and hence $Cl_P \simeq Cl_Q$. The Clifford algebra $Cl_Q$ is the $8$-dimensional $\C$-algebra generated by $x_1,x_2$ and $x_3$ satisfying the defining relations
\[
x_1^2=u,~x_2^2=v,~x_3^2=1 \qquad \text{and} \qquad x_ix_j+x_jx_i=0~\text{for $i \not= j$.} \]
If $uv \not= 0$ then $B_Q$ is a non-degenerate ternary quadratic form with determinant a square in $\C^*$ whence $Cl_Q$ is the direct sum of two copies of $M_2(\C)$. If $uv=0$, say $u=0$ and $v \not= 0$, then $x_1$ generates a nilpotent two-sided ideal of $Cl_Q$ and the quotient is the Clifford algebra of the non-degenerate binary quadratic form
\[
B_R = \begin{bmatrix} v & 0 \\ 0 & 1 \end{bmatrix} \qquad \text{whence} \qquad Cl_R \simeq M_2(\C) \]
as any such algebra is a quaternion algebra. Finally, if both $u=0=v$ then the two-sided ideal $I$ generated by $x_1$ and $x_2$ is nilpotent and the quotient
\[
Cl_R/I = \C[x_3]/(x_3^2-1) \simeq \C \oplus \C .\]
As the maximal ideals of a non-commutative algebra $R$ and of a quotient $R/I$ by a nilpotent ideal $I$ coincide, the statements follow.
\end{proof}

\begin{lemma} Intersecting with the center $R_c$ determines a continuous map
\[
\wis{spec}~\Lambda_c \rTo^{\psi} \wis{spec}~R_c, 
\]
which is a one-to-one correspondence away from the unique singularity of $\wis{spec}~R_c$ where the fiber consists of two points.
\end{lemma}

\begin{proof} The inclusion $\C[x,y,z] \subset R_c$ determines a two-fold cover
\[
\wis{spec}~R_c = \V(D^2-4(z^2-xy)) \subset \Af^4 \rOnto^c \Af^3 \qquad (x,y,z,D) \mapsto (x,y,z) \]
which is {\em ramified} over $\V(z^2-xy)$. That is, if $P=(a,b,c) \notin \V(z^2-xy)$ then there are exactly two points lying over it
\[
P_1 = (a,b,c,+\sqrt{c^2-ab}) \qquad \text{and} \qquad P_2 = (a,b,c,-\sqrt{c^2-ab}). \]
On the other hand, if $P = (a,b,c) \in \V(z^2-xy)$, then there is just one point lying over it : $(a,b,c,0)$. The statement then follows from combining this covering information with the composition map
\[
\wis{spec}~\Lambda_c \rTo^{\psi} \wis{spec}~R_c \rTo^c \Af^3 \]
which is $\phi$ in the foregoing lemma.
\end{proof}

Observe that $\psi$ is a homeomorphism on $\wis{spec}~\Lambda_c - \V(x,y,z)$ and hence can be seen as a non-commutative birational map. If $\mathfrak{m}$ lies in this open set then
\[
\Lambda_c/\mathfrak{m} \simeq M_2(\C) \]
whereas for the two maximal ideals $\mathfrak{m}_+ = (X,Y,Z-1)$ and $\mathfrak{m}_- = (X,Y,Z+1)$ lying over the conifold singularity we have
\[
\Lambda_c/\mathfrak{m}_+ \simeq \C \simeq \Lambda_c/\mathfrak{m}_-~.\]
We denote the associated one-dimensional $\Lambda_c$-representations by $\phi_+$ resp. $\phi_-$. It is easy to verify that these are the only one-dimensional $\Lambda_c$-representations.

\begin{proposition} For the conifold algebra $\Lambda_c$, the representation variety $\wis{rep}_2~\Lambda_c$ is a smooth affine variety having three disjoint irreducible components. Two of these components are a point, the third component $\wis{trep}_2~\Lambda_c$ has dimension $6$. In particular, the conifold algebra $\Lambda_c$ is a smooth order whence the birational map $\psi$ above can be viewed as a non-commutative desingularization.
\end{proposition}

\begin{proof} From the defining relation $Z^2 = 1$ it follows that the image of $Z$ in any finite dimensional representation has eigenvalues $\pm 1$. Hence, after simultaneous conjugation of the images of $X$, $Y$ and $Z$ we may assume that $Z$ has one of the following three forms
\[
Z \mapsto \begin{bmatrix} 1 & 0 \\ 0 & 1 \end{bmatrix} \quad \text{or} \quad Z \mapsto \begin{bmatrix} -1 & 0 \\ 0 & -1 \end{bmatrix} \quad \text{or} \quad Z \mapsto \begin{bmatrix} 1 & 0 \\ 0 & -1 \end{bmatrix}. \]
The first two possibilities are easily dealt with. Here, the image of $Z$ is a central unit so it follows from the relations $XZ+ZX=0=YZ+ZY$ as in the previous lemma that $X \mapsto 0$ and $Y \mapsto 0$. That is, these two components consist of just one point (the action of $\wis{GL}_2$ by simultaneous conjugation fixes these matrices) corresponding to the $2$-dimensional {\em semi-simple} representations
\[
M_+ = \phi_+ \oplus \phi_+ \qquad \text{and} \qquad M_- = \phi_- \oplus \phi_-~. \]
The interesting case is the third one. Because $X^2$ and $Y^2$ are central elements it follows (for example using the characteristic polynomial of $2 \times 2$ matrices) that in any $2$-dimensional representation $\Lambda_c \rTo^{\phi} M_2(\C)$ we have that $tr(\phi(X))=0$ and $tr(\phi(Y))=0$. Hence, the third component of $\wis{rep}_2~\Lambda_c$ consists of those $2$-dimensional representations $\phi$ such that
\[
tr(\phi(X)) = 0 \qquad tr(\phi(Y)) = 0 \qquad \text{and} \qquad tr(\phi(Z)) = 0. \]
For this reason we denote this component by $\wis{trep}_2~\Lambda_c$ and call it the variety of {\em trace preserving $2$-dimensional representations}. To describe the coordinate ring of this component we can use {\em trace zero} generic $2 \times 2$ matrices
\[
X \mapsto \begin{bmatrix} x_1 & x_2 \\ x_3 & -x_1 \end{bmatrix} \quad Y \mapsto \begin{bmatrix} y_1 & y_2 \\ y_3 & -y_1 \end{bmatrix} \quad Z \mapsto \begin{bmatrix} z_1 & z_2 \\ z_3 & -z_1 \end{bmatrix}
\]
which drastically reduces the defining equations as $T^2$ and $TS+ST$ are both scalar matrices for any trace zero $2 \times 2$ matrices. More precisely, we have
\[
XZ+ZX \mapsto \begin{bmatrix} 2x_1z_1+x_2z_3+x_3z_2 & 0 \\ 0 & 2x_1z_1+x_2z_3+x_3z_2 \end{bmatrix} \]
\[
YZ+ZY \mapsto \begin{bmatrix} 2y_1z_1+y_2z_3+y_3z_2 & 0 \\ 0 & 2y_1z_1+y_2z_3+y_3z_2 \end{bmatrix} \]
\[
Z^2 \mapsto \begin{bmatrix} z_1^2+z_2z_3 & 0 \\ 0 & z_1^2+z_2z_3 \end{bmatrix} \]
and therefore the coordinate ring of $\wis{trep}_2~\Lambda_c$
\[
\C[\wis{trep}_2~\Lambda_c] = \frac{\C[x_1,x_2,x_3,y_1,y_2,y_3,z_1,z_2,z_3]}{(2x_1z_1+x_2z_3+x_3z_2,2y_1z_1+y_2z_3+y_3z_2,z_1^2+z_2z_3-1)}. \]
To verify that $\wis{trep}_2~\Lambda_c$ is a smooth $6$-dimensional affine variety we therefore have to show that the {\em Jacobian matrix}
\[
\begin{bmatrix}
2z_1 & z_3 & z_2 & 0 & 0 & 0 & 2x_1 & x_3 & x_2 \\
0 & 0 & 0 & 2z_1 & z_3 & z_2 & 2y_1 & y_3 & y_2 \\
0 & 0 & 0 & 0 & 0 & 0 & 2z_1 & z_3 & z_2
\end{bmatrix}
\]
has constant rank $3$ on $\wis{trep}_2~\Lambda_c$. This is forced by the submatrices $\begin{bmatrix} 2z_1 & z_3 & z_2 \end{bmatrix}$ along the 'diagonal' of the Jacobian unless $z_1=z_2=z_3=0$ but this cannot hold for a point in $\wis{trep}_2~\Lambda_c$ by the equation
$z_1^2 + z_2z_3 = 1$.
\end{proof}

Next, we will use the two idempotents $e_1=\frac{1}{2}(Z-1)$ and $e_2=\frac{1}{2}(Z+1)$ to relate the conifold algebra $\Lambda_c$ to the quiver mentioned above. Consider a representation in $\wis{trep}_2~\Lambda_c$ then we can use base change to bring the image of $Z$ into the form
\[
Z \mapsto \begin{bmatrix} 1 & 0 \\ 0 & -1 \end{bmatrix} .\]
Taking the generic $2 \times 2$ matrices
\[
X \mapsto \begin{bmatrix} x_1 & x_2 \\ x_3 & x_4 \end{bmatrix} \qquad Y \mapsto \begin{bmatrix} y_1 & y_2 \\ y_3 & y_4 \end{bmatrix}
\]
it follows from the relations $XZ+ZX=0=YZ+ZY$ that $x_1=x_4=0=y_1=y_4$. Therefore, a representation in $\wis{trep}_2~\Lambda_c$ can be simultaneously conjugated to one of the form
\[
X \mapsto \begin{bmatrix} 0 & x_2 \\ x_3 & 0 \end{bmatrix} \quad Y \mapsto \begin{bmatrix}
0 & y_2 \\ y_3 & 0 \end{bmatrix} \quad Z \mapsto \begin{bmatrix} 1 & 0 \\ 0 & -1 \end{bmatrix} \]
and as the images of $X^2$ and $Y^2$ are scalar matrices the remaining defining relations $[X^2,Y]=0=[Y^2,X]$ are automatically satisfied. $2$-dimensional representations of $A_{con}$ in this canonical form hence form a smooth $4$-dimensional affine space
\[
\Af^4 = \V(x_1,x_4,y_1,y_4,z_1-1,z_2,z_3,z_4+1) \subset \Af^{12}. \]
To recover $\wis{trep}_2~\Lambda_c$ from this affine space we let $\wis{GL}_2$ act on it. The subgroup of $\wis{GL}_2$ fixing the matrix
\[
\begin{bmatrix} 1 & 0 \\ 0 & -1 \end{bmatrix} \qquad \text{is} \qquad T = \{ \begin{bmatrix} \lambda & 0 \\ 0 & \mu \end{bmatrix}~|~\lambda,\mu \in \C^* \}, \]
the two-dimensional  torus. There is an action of $T$ on the product $\wis{GL}_2 \times \Af^4$ via
\[
t.(g,P) = (gt^{-1},t.P) \qquad \text{for all $t \in T, g \in \wis{GL}_2$ and $P \in \Af^4$} \]
and where $t.P$ means the action by simultaneous conjugation by the $2 \times 2$ matrix $t \in T \subset \wis{GL}_2$ on the three $2 \times 2$ matrix-components of $P$.

\begin{proposition} Under the action-map 
\[
\wis{GL}_2 \times \Af^4 \rTo \wis{trep}_2~\Lambda_c \qquad (g,P) \mapsto g.P \]
two points $(g,P)$ and $(g',P')$ are mapped to the same point if and only if they belong to the same $T$-orbit in $\wis{GL}_2 \times \Af^4$. That is, we can identify $\wis{trep}_2~\Lambda_{c}$ with the principal fiber bundle
\[
\wis{trep}_2~\Lambda_c \simeq \wis{GL}_2 \times^T \Af^4 = (\wis{GL}_2 \times \Af^4) / T. \]
In particular, there is a natural one-to-one correspondence between $\wis{GL}_2$-orbits in $\wis{trep}_2~\Lambda_c$ and $T$-orbits in $\Af^4$. Observe that one can identify the $T$-action on $\Af^4$ with the $\wis{GL}(\alpha)$-action on the representation space $\wis{rep}_{\alpha}~Q$ for the quiver-setting
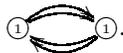
\[
\xymatrix{\vtx{1} \ar@/^/[r] \ar@/^2ex/[r] & \vtx{1} \ar@/^/[l] \ar@/^2ex/[l]}.
\]
In particular, the conifold algebra $\Lambda_c$ is the quiver-order $\int_{\alpha}~\C Q$.
\end{proposition}

\begin{proof} If $g.P = g'.P'$, then $P = g^{-1}g'.P'$ and as both $P$ and $P'$ have as their third $2 \times 2$ matrix component 
\[
\begin{bmatrix} 1 & 0 \\ 0 & -1 \end{bmatrix} \]
it follows that $g^{-1}g'$ is in the stabilizer subgroup of this matrix so $g^{-1}g' = t^{-1}$ for some $t \in T$
whence $g' = gt^{-1}$ and as $(g^{-1}g')^{-1}.P = P'$ also $t.P = P'$ whence
\[
t.(g,P) = (gt^{-1},t.P) = (g',P') \]
Hence we can identify $\wis{trep}_2~\Lambda_c = \wis{GL}_2.\Af^4$ with the orbit-space of the $T$-action which is just $\wis{GL}_2 \times^T \Af^4$. Incidentally, this gives another proof for smoothness of $\wis{trep}_2~\Lambda_c$ as it is the base of a fibration with smooth fibers of the smooth top space $\wis{GL}_2 \times \Af^4$.
$\wis{GL}_2$ acts on $\wis{GL}_2 \times \Af^4$ by $g.(g',P') = (gg',P')$ and this action commutes with the $T$-action so induces a $\wis{GL}_2$-action on the orbit-space
\[
\wis{GL}_2 \times (\wis{GL}_2 \times^T \Af^4) \rTo \wis{GL}_2 \times^T \Af^4 \qquad g.\overline{(g',P')} = \overline{(gg',P')}.
\]
As we have identified $\wis{GL}_2 \times^T \Af^4$ with $\wis{trep}_2~\Lambda_c$ via the action map, that is
$\overline{(g,P)} = g.P$ the remaining statements follow.
\end{proof}

In this specific case we can explicitly compute polynomial (semi)-invariants using the $T$-action and relate it to the general results mentioned before.

\begin{lemma} The ring of polynomial invariants
\[
\C[\wis{trep}_2~\Lambda_c]^{\wis{GL}_2} \simeq \C[\Af^4]^T \]
is isomorphic to the coordinate ring of the conifold singularity $R_c$ and the quotient map
\[
\wis{trep}_2~\Lambda_c \rOnto \wis{spec}~R_c \]
maps a two-dimensional representation to the direct sum of its Jordan-H\"older components. \end{lemma}

\begin{proof} The action of the two-dimensional torus $T$ on $\Af^4 = \{ (x_2,x_3,y_2,y_3) \}$ is given by
\[
\begin{bmatrix} \lambda & 0 \\ 0 & \mu \end{bmatrix}.(\begin{bmatrix} 0 & x_2 \\ x_3 & 0 \end{bmatrix},
\begin{bmatrix} 0 & y_2 \\ y_3 & 0 \end{bmatrix},\begin{bmatrix} 1 & 0 \\ 0 & -1 \end{bmatrix} ) = \]
\[
(\begin{bmatrix} 0 &  \lambda\mu^{-1}x_2 \\ \lambda^{-1}\mu x_3 & 0 \end{bmatrix},
\begin{bmatrix} 0 & \lambda \mu^{-1}y_2 \\ \lambda^{-1} \mu y_3 & 0 \end{bmatrix},\begin{bmatrix} 1 & 0 \\ 0 & -1 \end{bmatrix} ) .
\]
Hence, the action of $(\lambda,\mu) \in T$  on $\C[\Af^4] = \C[X_2,X_3,Y_2,Y_3]$ is defined by
\[
X_2 \mapsto \lambda^{-1}\mu X_2 \quad X_3 \mapsto \lambda \mu^{-1} X_3 \quad Y_2 \mapsto \lambda^{-1}\mu Y_2 \quad Y_3 \mapsto \lambda \mu^{-1} Y_3 \]
and this action sends any monomial in the variables to a scalar multiple of that monomial. So, in order to determine the ring of polynomial invariants
\[
\C[X_2,X_3,Y_2,Y_3]^T = \{ f= f(X_2,X_3,Y_2,Y_3)~|~(\lambda,\mu).f = f~ \forall (\lambda,\mu) \in T \}
\]
it suffices to determine all invariant monomials, or equivalently, all positive integer quadruplets $(a,b,c,d)$ such that $a-b+c-d=0$ as
\[
(\lambda,\mu).X_2^aX_3^bY_2^cY_3^d = \lambda^{-a+b-c+d} \mu^{a-b+c-d} X_2^aX_3^bY_2^cY_3^d
\]
Clearly, such quadruplets are all generated (as Abelian group under addition) by the four basic ones
\[
(1,1,0,0) \mapsto X_2X_3 \quad (1,0,0,1) \mapsto X_2Y_3 \quad (0,1,1,0) \mapsto X_3Y_2 \quad (0,0,1,1) \mapsto Y_2Y_3 \]
and therefore
\[
\C[\wis{trep}_2~\Lambda_c]^{\wis{GL}_2} \simeq \C[X_2,X_3,Y_2,Y_3]^T = \C[X_2X_3,X_2Y_3,X_3Y_2,Y_2,Y_3] \simeq \frac{\C[p,q,r,s]}{(ps-qr)} \]
is the conifold singularity $R_c$. We know already that $\wis{spec}~R_c$ has as its points the isomorphism classes of $2$-dimensional semi-simple representations with $\phi_+ \oplus \phi_-$ as the semi-simple representation corresponding to the singularity and all other points classify a unique simple $2$-dimensional representation. 
\end{proof}

For the quiver-setting $(Q,\alpha)$ there are essentially two stability structures : $\theta=(-1,1)$ and $\theta'=(1,-1)$. Again, we can use elementary arguments in this case to calculate the moduli spaces.

\begin{lemma} The moduli space of all $\theta$-(semi)stable $\alpha$-dimensional representations
\[
\wis{moduli}^{\theta}_{\alpha}~\Lambda_c \simeq \wis{proj}~\C[\wis{rep}_{\alpha}~Q]^{\wis{GL}(\alpha),\theta} \]
is the $\wis{proj}$ of the ring of $\theta$-semi-invariants and as the semi-invariants of weight zero are the polynomial invariants we get a projective morphism
\[
\wis{proj}~\C[\wis{rep}_{\alpha}~Q]^{\wis{GL}(\alpha),\theta} \rOnto \wis{spec}~R_c \]
which is a desingularization of the conifold singularity.
\end{lemma}

\begin{proof}
As in the case of polynomial invariants, the space $\C[\wis{rep}_{\alpha}~Q]^{\wis{GL}(\alpha),\theta}_k$ is spanned by monomials
\[
x_2^ax_3^by_2^cy_3^d \qquad \text{satisfying} \qquad -a+b-c+d=k \]
and one verifies that this space is the module over the ring of polynomial invariants generated by all monomials of degree $k$ in $x_3$ and $y_3$. That is
\[
\C[\wis{rep}_{\alpha}~Q]^{\wis{GL}(\alpha),\theta} = \C[x_2x_3,x_2y_3,x_3y_2,y_2y_3][x_3,y_3] \subset \C[x_2,y_2,x_3,y_3] \]
with the generators $a=x_2x_3,b=x_2y_3,c=x_3y_2$ and $d=y_2y_3$ of degree zero and $e=x_3$ and $f=y_3$ of degree one. As a consequence, we can identify $\wis{proj}~\C[\wis{rep}_{\alpha}~Q]^{\wis{GL}(\alpha),\theta}$ with the closed subvariety
\[
\V(ad-bc,af-be,cf-de) \subset \Af^4 \times \PP^1 \]
with $(a,b,c,d)$ the affine coordinates of $\Af^4$ and $[e:f]$ projective coordinates of $\PP^1$. The projection $\wis{proj}~\C[\wis{rep}_{\alpha}~Q]^{\wis{GL}(\alpha),\theta} \rOnto \wis{spec}~R_c$ is projection onto the $\Af^4$-component.

To prove smoothness we cover $\PP^1$ with the two affine opens $e \not= 0$ (with affine coordinate $x = f/e$ and $f \not=0$ with affine coordinate $y = e/f$. In the affine coordinates $(a,b,c,d,x)$ the relations become
\[
ad = bc \qquad ax = b \qquad \text{and} \qquad cx = d \]
whence the coordinate ring is $\C[a,c,x]$ and so the variety is smooth on this affine open. Similarly, the coordinate ring on the other affine open is $\C[b,d,y]$ and smoothness follows. Moreover, $\pi$ is {\em birational} over the complement of the singularity. This follows from the relations
\[
ax = b, \quad cx = d, \quad by = a, \quad dy = c \]
which determine $x$ (or $y$ and hence the point in $\wis{proj}$) lying over any $(a,b,c,d) \not= (0,0,0,0)$ in $\wis{spec}~R_c$. Therefore, the map $\pi$ is a desingularization and the {\em exceptional fiber}
\[
E = \pi^{-1}(0,0,0,0) \simeq \PP^1 \]
which classifies the $\theta$-stable representations which lie over $(0,0,0,0)$ (that is, those such that $x_2x_3=x_2y_3=x_3y_2=y_2y_3=0$) as they are all of the form
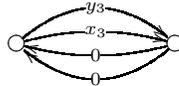
\[
\xymatrix{\vtx{} \ar@/^1ex/[rr]|{x_3} \ar@/^3ex/[rr]|{y_3} & & \vtx{} \ar@/^1ex/[ll]|{0} \ar@/^3ex/[ll]|{0}}
\]
with either $x_3 \not= 0$ or $y_3 \not= 0$ and the different $T$-orbits of those are parametrized by the points of $\PP^1$. As the smooth points of $\wis{spec}~R_c$ are known to correspond to isomorphism classes of simple (hence certainly $\theta$-stable) representations we have proved that
\[
\wis{proj}~\C[\wis{rep}_{\alpha}~Q]^{\wis{GL}(\alpha),\theta} \simeq \wis{moduli}^{\theta}_{\alpha}~\Lambda_c \]
is the moduli space of all $\theta$-stable $\alpha$-dimensional representations of $Q$.
\end{proof}

Clearly, we could have done the same calculations starting with the stability structure $\theta' = (1,-1)$ and obtained another desingularization replacing the roles of $x_2,y_2$ and $x_3,y_3$. This gives us the situation
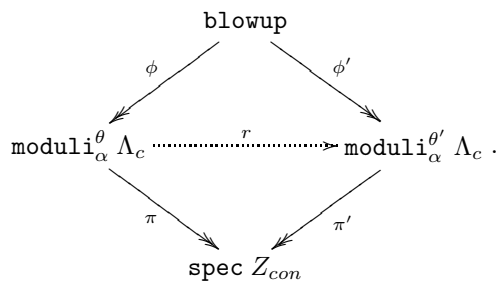
\[
\vcenter{\xymatrix@C=10pt@R=30pt{
& \wis{blowup} \ar@{->>}[dl]_{\phi} \ar@{->>}[dr]^{\phi'}\\
\wis{moduli}^{\theta}_{\alpha}~\Lambda_c \ar@{->>}[dr]_{\pi} \ar@{.>}[rr]^r& & \wis{moduli}^{\theta'}_{\alpha}~\Lambda_c \ar@{->>}[dl]^{\pi'}\\
& \wis{spec}~Z_{con}
}}.
\]
Here, $\wis{blowup}$ denotes the desingularization of $\wis{spec}~R_c$ one obtains by blowing-up the point $(0,0,0,0) \in \Af^4$ and which has exceptional fiber $\PP^1 \times \PP^1$. Blowing down either of these lines (the maps $\phi$ and $\phi'$) one obtains the 'minimal' resolutions given by the moduli spaces. These spaces are related by the {\em rational map} $r$ which is called the {\em Atiyah flop} in string theory-literature.

 \end{document}